\documentclass[11pt]{amsart}

\usepackage{amssymb}
\usepackage{graphicx}
\usepackage[all]{xy}

\title[Invariance of the branching and merging homologies]{T-homotopy and refinement of observation (III) : Invariance of the branching and merging homologies}
\author[P. Gaucher]{Philippe Gaucher}
\address{Preuves Programmes et Syst{\`e}mes\\ Universit{\'e} Paris 7--Denis Diderot\\
Case 7014\\2 Place Jussieu\\ 75251 PARIS Cedex 05\\ France}
\email{gaucher@pps.jussieu.fr \\ 
http://www.pps.jussieu.fr/{\~{}}gaucher/}
\subjclass{55U35, 55P99, 68Q85} 
\keywords{concurrency, homotopy, directed homotopy, model category, refinement of observation, poset, cofibration, Reedy category, homotopy colimit, branching, merging, homology}

\newcommand{\C}{\mathcal{C}}

\newcommand{\Z}{\mathbb{Z}}
\newcommand{\de}{\partial}
\newcommand{\p}\times
\renewcommand{\vec}{\overrightarrow}
\renewcommand{\P}{\mathbb{P}}

\newcommand{\be}{\begin{equation}}
\newcommand{\ee}{\end{equation}}
\newcommand{\bea}{\begin{eqnarray}}
\newcommand{\eea}{\end{eqnarray}}
\newcommand{\beas}{\begin{eqnarray*}}
\newcommand{\eeas}{\end{eqnarray*}}

\newtheorem*{thmN}{Theorem}

\newtheorem{thm}{Theorem}[section]
\newtheorem{prop}[thm]{Proposition}
\newtheorem{lem}[thm]{Lemma}
\newtheorem{cor}[thm]{Corollary}

\newtheorem{defn}[thm]{Definition}
\newtheorem{nota}[thm]{Notation}
\newcommand{\bd}{\begin{defn}}
\newcommand{\ed}{\end{defn}}
\newcommand{\bp}{\begin{prop}}
\newcommand{\ep}{\end{prop}}
\newcommand{\bth}{\begin{thm}}
\renewcommand{\eth}{\end{thm}}
\newcommand{\bpf}{\begin{proof}}
\newcommand{\epf}{\end{proof}}

\newcommand{\fl}[1]{\ar@{->}[l]_{#1}}
\newcommand{\fr}[1]{\ar@{->}[r]^-{#1}}
\newcommand{\fd}[1]{\ar@{->}[d]_{#1}}
\newcommand{\fu}[1]{\ar@{->}[u]^{#1}}
\newcommand{\f}[2]{\ar@{->}[#1]|{#2}}
\newcommand{\ff}[2]{\ar@2{->}[#1]|{#2}}
\newcommand{\frr}[1]{\ar@{->}[rr]^{#1}}

\newcommand{\iso}{\cong}
\newcommand{\lp}{\left(}
\newcommand{\rp}{\right)}

\newcommand{\vI}{\vec{I}}
\renewcommand{\leq}{\leqslant}
\renewcommand{\geq}{\geqslant}

\def\cartesien{%
  \ar@{-}[]+R+<6pt,-2pt>;[]+RD+<6pt,-6pt>%
  \ar@{-}[]+D+<2pt,-6pt>;[]+RD+<6pt,-6pt>%
}
\def\cocartesien{%
  \ar@{-}[]+L+<-6pt,+2pt>;[]+LU+<-6pt,+6pt>%
  \ar@{-}[]+U+<-2pt,+6pt>;[]+LU+<-6pt,+6pt>%
}

\newcommand{\brm}[1]{\rm{\mathbf{#1}}}
\renewcommand{\top}{{\brm{Top}}}
\newcommand{\dtop}{{\brm{Flow}}}
\newcommand{\set}{{\brm{Set}}}
\newcommand{\glob}{{\rm{Glob}}}
\DeclareMathOperator{\sing}{Sing}
\DeclareMathOperator{\obj}{Obj}

\newcommand{\liminj}{\varinjlim}

\makeatletter

\def\varholim@#1#2{%
  \vtop{\m@th\ialign{##\cr
    \hfil$#1\operator@font holim$\hfil\cr
    \noalign{\nointerlineskip\kern1.5\ex@}#2\cr
    \noalign{\nointerlineskip\kern-\ex@}\cr}}%
}
\def\holimproj{%
  \mathop{\mathpalette\varholim@{\leftarrowfill@\textstyle}}\nmlimits@
}
\def\holiminj{%
  \mathop{\mathpalette\varholim@{\rightarrowfill@\textstyle}}\nmlimits@
}

\makeatother

\DeclareMathOperator{\hop}{ho\P} 
\DeclareMathOperator{\id}{Id}
\DeclareMathOperator{\card}{card}

\DeclareMathOperator{\diag}{{\rm{Diag}}}

\DeclareMathOperator{\cell}{{\brm{cell}}}
\DeclareMathOperator{\cof}{{\brm{cof}}}
\DeclareMathOperator{\inj}{{\brm{inj}}}
\newcommand{\hda}{{\cell(\dtop)}}
\DeclareMathOperator{\im}{im}

\begin{document}

\begin{abstract} 
  This series explores a new notion of T-homotopy equivalence of
  flows.  The new definition involves embeddings of finite bounded
  posets preserving the bottom and the top elements and the associated
  cofibrations of flows. In this third part, it is proved that the
  generalized T-homotopy equivalences preserve the bran\-ching and
  merging homology theories of a flow. These homology theories are of
  interest in computer science since they detect the non-deterministic
  branching and merging areas of execution paths in the time flow of a
  higher dimensional automaton. The proof is based on Reedy model
  category techniques.
\end{abstract}

\maketitle

\tableofcontents

\section{Outline of the paper}

The main feature of the algebraic topological model of \textit{higher
  dimensional automata} (or HDA) introduced in \cite{model3}, the
category of \textit{flows}, is to provide a framework for modelling
continuous deformations of HDA corresponding to subdivision or
refinement of observation.  The equivalence relation induced by these
deformations, called \textit{dihomotopy}, preserves geometric
properties like the \textit{initial} or \textit{final states}, and
therefore computer-scientific properties like the presence or not of
\textit{deadlocks} or of \textit{unreachable states} in concurrent
systems \cite{survol}.  More generally, dihomotopy is designed to
preserve all computer-scientific property invariant by refinement of
observation.  Figure~\ref{ex1} represents a very simple example of
refinement of observation, where a $1$-dimensional transition from an
initial state to a final state is identified with the composition of
two such transitions.

In the framework of flows, there are two kinds of dihomotopy
equivalences \cite{ConcuToAlgTopo}: the \textit{weak S-homotopy
  equivalences} (the spatial deformations of \cite{ConcuToAlgTopo})
and the \textit{T-homotopy equivalences} (the temporal deformations of
\cite{ConcuToAlgTopo}). The geometric explanations underlying the
intuition of S-homotopy and T-homotopy are given in the first part of
this series \cite{1eme}, but the reference \cite{diCW} must be
preferred.

It is very fortunate that the class of weak S-homotopy equivalences
can be interpreted as the class of weak equivalences of a model
structure \cite{model3} in the sense of Hovey's book
\cite{MR99h:55031}.  This fact makes their study easier.  Moreover,
this model structure is necessary for the formulation of the only
known definition of T-homotopy.

The purpose of this paper is to prove that the new notion of
T-homotopy equivalence is well-behaved with respect to the
\textit{branching and merging homologies} of a flow. The latter
homology theories are able to detect the non-deterministic higher
dimensional branching and merging areas of execution paths in the time
flow of a higher dimensional automaton \cite{exbranch}. More
precisely, one has:
\begin{thmN} (Corollary~\ref{enfinfin})
Let $f:X\longrightarrow Y$ be a generalized T-homotopy equivalence.
Then for any $n\geq 0$, the morphisms of abelian groups
$H_n^-(f):H_n^-(X)\longrightarrow H_n^-(Y)$,
$H_n^+(f):H_n^+(X)\longrightarrow H_n^+(Y)$ are isomorphisms of groups
where $H_n^-$ (resp. $H_n^+$) is the $n$-th branching (resp. merging)
homology group.
\end{thmN}
The core of the paper starts with Section~\ref{remflow} which recalls
the definition of a flow and the description of the weak S-homotopy
model structure. The latter is a fundamental tool for the
sequel. Section~\ref{section_defT} recalls the new notion of
T-homotopy equivalence.

Section~\ref{principlemainproof} recalls the definition of the
branching space and the homotopy branching space of a flow. The same
section explains the principle of the proof of the following theorem:
\begin{thmN} (Theorem~\ref{resultat1})
The homotopy branching space of a full directed ball at any state
different from the final state is contractible (it is empty at the
final state).
\end{thmN}
We give the idea of the proof for a full directed ball which is not
too simple, and not too complicated.  The latter theorem is the
technical core of the paper because a generalized T-homotopy
equivalence consists in replacing in a flow a full directed ball by a
more refined full directed ball (Figure~\ref{ex2}), and in iterating
this replacement process transfinitely.

Section~\ref{maintool} introduces a diagram of topological spaces
$\mathcal{P}_\alpha^-(X)$ whose colimit calculates the branching space
$\P_\alpha^- X$ for every loopless flow $X$ (Theorem~\ref{calculPm})
and every $\alpha\in X^0$.  Section~\ref{reedythree} builds a Reedy
structure on the base category of the diagram
$\mathcal{P}_\alpha^-(X)$ for any loopless flow $X$ whose poset
$(X^0,\leq)$ is locally finite so that the colimit functor becomes a
left Quillen functor (Theorem~\ref{equi}). Section~\ref{diagcofibrant}
then shows that the diagram $\mathcal{P}_\alpha^-(X)$ is Reedy
cofibrant as soon as $X$ is a cell complex of the model category
$\dtop$ (Theorem~\ref{cascofibrant}).  Section~\ref{end} completes the
proof that the homotopy branching and homotopy merging spaces of every
full directed ball are contractible (Theorem~\ref{resultat1}).
Section~\ref{rappelHMP} recalls the definition of the branching and
merging homology theories. At last, Section~\ref{preHMP} proves the
invariance of the branching and merging homology theories with respect
to T-homotopy.

\vspace{1cm} \textit{Warning.} This paper is the third part of a
series of papers devoted to the study of T-homotopy. Several other
papers explain the geometrical content of T-homotopy. The best
reference is probably \cite{diCW} (it does not belong to the series).
The knowledge of the first and second parts is not required, except
for the left properness of the weak S-homotopy model structure of
$\dtop$ available in \cite{2eme}. The latter fact is used twice in the
proof of Theorem~\ref{pre1}. The material collected in the appendices
\ref{limelm}, \ref{calpush} and \ref{mixcomp} will be reused in the
fourth part \cite{4eme}. The proofs of these appendices are
independent from the technical core of this part.

\section{Prerequisites and notations}

The initial object (resp. the terminal object) of a category $\C$, if
it exists, is denoted by $\varnothing$ (resp. $\mathbf{1}$).

Let $\C$ be a cocomplete category.  If $K$ is a set of morphisms of
$\C$, then the class of morphisms of $\C$ that satisfy the RLP
(\textit{right lifting property}) with respect to any morphism of $K$
is denoted by $\inj(K)$ and the class of morphisms of $\C$ that are
transfinite compositions of pushouts of elements of $K$ is denoted by
$\cell(K)$. Denote by $\cof(K)$ the class of morphisms of $\C$ that
satisfy the LLP (\textit{left lifting property}) with respect to the
morphisms of $\inj(K)$.  It is a purely categorical fact that
$\cell(K)\subset \cof(K)$. Moreover, every morphism of $\cof(K)$ is a
retract of a morphism of $\cell(K)$ as soon as the domains of $K$ are
small relative to $\cell(K)$ (\cite{MR99h:55031} Corollary~2.1.15). An
element of $\cell(K)$ is called a \textit{relative $K$-cell complex}.
If $X$ is an object of $\C$, and if the canonical morphism
$\varnothing\longrightarrow X$ is a relative $K$-cell complex, then
the object $X$ is called a \textit{$K$-cell complex}.

Let $\C$ be a cocomplete category with a distinguished set of
morphisms $I$. Then let $\cell(\C,I)$ be the full subcategory of $\C$
consisting of the object $X$ of $\C$ such that the canonical morphism
$\varnothing\longrightarrow X$ is an object of $\cell(I)$. In other
terms, $\cell(\C,I)=(\varnothing\!\downarrow \! \C) \cap \cell(I)$.

It is obviously impossible to read this paper without a strong
familiarity with \textit{model categories}. Possible references for
model categories are \cite{MR99h:55031}, \cite{ref_model2} and
\cite{MR1361887}.  The original reference is \cite{MR36:6480} but
Quillen's axiomatization is not used in this paper. The axiomatization
from Hovey's book is preferred.  If $\mathcal{M}$ is a
\textit{cofibrantly generated} model category with set of generating
cofibrations $I$, let $\cell(\mathcal{M}) := \cell(\mathcal{M},I)$ :
this is the full subcategory of \textit{cell complexes} of the model
category $\mathcal{M}$. A cofibrantly generated model structure
$\mathcal{M}$ comes with a \textit{cofibrant replacement functor}
$Q:\mathcal{M} \longrightarrow \cell(\mathcal{M})$. For any morphism
$f$ of $\mathcal{M}$, the morphism $Q(f)$ is a cofibration, and even
an inclusion of subcomplexes (\cite{ref_model2} Definition~10.6.7)
because the cofibrant replacement functor $Q$ is obtained by the small
object argument.

A \textit{partially ordered set} $(P,\leq)$ (or \textit{poset}) is a
set equipped with a reflexive antisymmetric and transitive binary
relation $\leq$. A poset is \textit{locally finite} if for any
$(x,y)\in P\p P$, the set $[x,y]=\{z\in P,x\leq z\leq y\}$ is finite.
A poset $(P,\leq)$ is \textit{bounded} if there exist $\widehat{0}\in
P$ and $\widehat{1}\in P$ such that $P = [\widehat{0},\widehat{1}]$
and such that $\widehat{0} \neq \widehat{1}$. Let $\widehat{0}=\min P$
(the bottom element) and $\widehat{1}=\max P$ (the top element). In a
poset $P$, the interval $]\alpha,-]$ (the sub-poset of elements of $P$
strictly bigger than $\alpha$) can also be denoted by $P_{>\alpha}$.

A poset $P$, and in particular an ordinal, can be viewed as a small
category denoted in the same way: the objects are the elements of $P$
and there exists a morphism from $x$ to $y$ if and only if $x\leq
y$. If $\lambda$ is an ordinal, a \textit{$\lambda$-sequence} in a
cocomplete category $\C$ is a colimit-preserving functor $X$ from
$\lambda$ to $\C$. We denote by $X_\lambda$ the colimit $\liminj X$
and the morphism $X_0\longrightarrow X_\lambda$ is called the
\textit{transfinite composition} of the $X_\mu\longrightarrow
X_{\mu+1}$.

Let $\C$ be a category. Let $\alpha$ be an object of $\C$. The
\textit{latching category} $\de(\C\!\downarrow\! \alpha)$ at $\alpha$ is the
full subcategory of $\C\!\downarrow\! \alpha$ containing all the
objects except the identity map of $\alpha$. The
\textit{matching category} $\de(\alpha\!\downarrow\!\C)$ at $\alpha$ is the
full subcategory of $\alpha\!\downarrow\!\C$ containing all the
objects except the identity map of $\alpha$.

Let $\mathcal{B}$ be a small category. A \textit{Reedy structure} on
$\mathcal{B}$ consists of two subcategories $\mathcal{B}_-$ and
$\mathcal{B}_+$, a map $d:\obj(\mathcal{B})\longrightarrow \lambda$
from the set of objets of $\mathcal{B}$ to some ordinal $\lambda$
called the {\rm degree function}, such that every non identity map in
$\mathcal{B}_+$ raises the degree, every non identity map in
$\mathcal{B}_-$ lowers the degree, and every map $f\in \mathcal{B}$
can be factored uniquely as $f=g \circ h$ with $h \in \mathcal{B}_-$
and $g \in \mathcal{B}_+$. A small category together with a Reedy
structure is called a \textit{Reedy category}.

If $\C$ is a small category and if $\mathcal{M}$ is a category, the
notation $\mathcal{M}^\C$ is the category of functors from $\C$ to
$\mathcal{M}$, i.e. the category of diagrams of objects of
$\mathcal{M}$ over the small category $\C$.

Let $\C$ be a complete and cocomplete category. Let $\mathcal{B}$ be a
Reedy category. Let $i$ be an object of $\mathcal{B}$. The \textit{latching
space functor} is the composite $L_i:\C^\mathcal{B}\longrightarrow
\C^{\de(\mathcal{B}_+\!\downarrow\! i)}\longrightarrow \C$ where the latter
functor is the colimit functor.  The \textit{matching space functor}
is the composite $M_i:\C^\mathcal{B}\longrightarrow
\C^{\de(i\!\downarrow\!\mathcal{B}_-)}\longrightarrow \C$ where the latter
functor is the limit functor.

A model category is \textit{left proper} if the pushout of a weak
equivalence along a cofibration is a weak equivalence. The model
categories $\top$ and $\dtop$ (see below) are both left proper.

In this paper, the notation $\xymatrix@1{\ar@{^{(}->}[r]&}$ means
\textit{cofibration}, the notation $\xymatrix@1{\ar@{->>}[r]&}$ means
\textit{fibration}, the notation $\simeq$ means \textit{weak
  equivalence}, and the notation $\iso$ means \textit{isomorphism}.

A categorical adjunction $\mathbb{L}:\mathcal{M}\leftrightarrows
\mathcal{N}:\mathbb{R}$ between two model categories is a
\textit{Quillen adjunction} if one of the following equivalent
conditions is satisfied: 1) $\mathbb{L}$ preserves cofibrations and
trivial cofibrations, 2) $\mathbb{R}$ preserves fibrations and trivial
fibrations. In that case, $\mathbb{L}$ (resp. $\mathbb{R}$) preserves
weak equivalences between cofibrant (resp. fibrant) objects.

If $P$ is a poset, let us denote by $\Delta(P)$ the \textit{order
  complex} associated with $P$. Recall that the order complex is a
simplicial complex having $P$ as underlying set and having the subsets
$\{x_0,x_1,\dots,x_n\}$ with $x_0<x_1<\dots<x_n$ as $n$-simplices
\cite{MR493916}.  Such a simplex will be denoted by
$(x_0,x_1,\dots,x_n)$.  The order complex $\Delta(P)$ can be viewed as
a poset ordered by the inclusion, and therefore as a small category.
The corresponding category will be denoted in the same way. The
opposite category $\Delta(P)^{op}$ is freely generated by the
morphisms $\de_i:(x_0,\dots,x_n) \longrightarrow
(x_0,\dots,\widehat{x_i},\dots,x_n)$ for $0\leq i\leq n$ and by the
simplicial relations $\de_i\de_j=\de_{j-1}\de_i$ for any $i<j$, where
the notation $\widehat{x_i}$ means that $x_i$ is removed.

If $\C$ is a small category, then the \textit{classifying space} of $\C$ is
denoted by $B\C$ \cite{MR0232393} \cite{MR0338129}.

The category $\top$ of \textit{compactly generated topological spaces}
(i.e. of weak Hausdorff $k$-spaces) is complete, cocomplete and
cartesian closed (more details for this kind of topological spaces in
\cite{MR90k:54001,MR2000h:55002}, the appendix of \cite{Ref_wH} and
also the preliminaries of \cite{model3}). For the sequel, all 
topological spaces will be supposed to be compactly generated. A
\textit{compact space} is always Hausdorff.

\section{Reminder about the category of flows}
\label{remflow}

The category $\top$ is equipped with the unique model structure having
the \textit{weak homotopy equivalences} as weak equivalences and having the
\textit{Serre fibrations}~\footnote{that is a continuous map having the RLP
with respect to the inclusion $\mathbf{D}^n\p 0\subset
\mathbf{D}^n\p [0,1]$ for any $n\geq 0$ where $\mathbf{D}^n$ is the
$n$-dimensional disk.} as fibrations.

The time flow of a higher dimensional automaton is encoded in an
object called a \textit{flow} \cite{model3}. A flow $X$ consists of a
set $X^0$ called the \textit{$0$-skeleton} and whose elements
correspond to the \textit{states} (or \textit{constant execution
  paths}) of the higher dimensional automaton. For each pair of states
$(\alpha,\beta)\in X^0\p X^0$, there is a topological space
$\P_{\alpha,\beta}X$ whose elements correspond to the
\textit{(non-constant) execution paths} of the higher dimensional
automaton \textit{beginning at} $\alpha$ and \textit{ending at}
$\beta$. For $x\in \P_{\alpha,\beta}X$, let $\alpha=s(x)$ and
$\beta=t(x)$. For each triple $(\alpha,\beta,\gamma)\in X^0\p X^0\p
X^0$, there exists a continuous map $*:\P_{\alpha,\beta}X\p
\P_{\beta,\gamma}X\longrightarrow \P_{\alpha,\gamma}X$ called the
\textit{composition law} which is supposed to be associative in an
obvious sense. The topological space $\P
X=\bigsqcup_{(\alpha,\beta)\in X^0\p X^0}\P_{\alpha,\beta}X$ is called
the \textit{path space} of $X$. The category of flows is denoted by
$\dtop$. A point $\alpha$ of $X^0$ such that there are no non-constant
execution paths ending at $\alpha$ (resp. starting from $\alpha$) is
called an \textit{initial state} (resp. a \textit{final state}). A
morphism of flows $f$ from $X$ to $Y$ consists of a set map $f^0:X^0
\longrightarrow Y^0$ and a continuous map $\P f: \P X \longrightarrow
\P Y$ preserving the structure. A flow is therefore ``almost'' a small
category enriched in $\top$.

\begin{figure}
\begin{center}
\includegraphics[width=7cm]{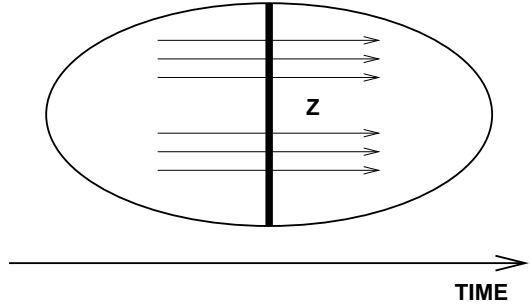}
\end{center}
\caption{Symbolic representation of
$\glob(Z)$ for some topological space $Z$} \label{exglob}
\end{figure}

An important example is the flow $\glob(Z)$ defined by \beas
&& \glob(Z)^0=\{\widehat{0},\widehat{1}\} \\
&& \P \glob(Z)=Z \\
&& s=\widehat{0} \\
&& t=\widehat{1} \eeas and a trivial composition law (cf.
Figure~\ref{exglob}).

The category $\dtop$ is equipped with the unique model structure
such that \cite{model3}: 
\begin{itemize}
\item The weak equivalences are the \textit{weak S-homotopy equivalences}, 
i.e. the morphisms of flows $f:X\longrightarrow Y$ such that
$f^0:X^0\longrightarrow Y^0$ is a bijection and such that $\P f:\P
X\longrightarrow \P Y$ is a weak homotopy equivalence. 
\item The fibrations are the morphisms of flows
$f:X\longrightarrow Y$ such that $\P f:\P X\longrightarrow \P Y$ is a
Serre fibration. 
\end{itemize}
This model structure is cofibrantly generated. The set of generating
cofibrations is the set $I^{gl}_+=I^{gl}\cup
\{R:\{0,1\}\longrightarrow \{0\},C:\varnothing\longrightarrow \{0\}\}$
with
\[I^{gl}=\{\glob(\mathbf{S}^{n-1})\subset \glob(\mathbf{D}^{n}), n\geq
0\}\] where $\mathbf{D}^{n}$ is the $n$-dimensional disk and
$\mathbf{S}^{n-1}$ the $(n-1)$-dimensional sphere. The set of
generating trivial cofibrations is
\[J^{gl}=\{\glob(\mathbf{D}^{n}\p\{0\})\subset
\glob(\mathbf{D}^{n}\p [0,1]), n\geq 0\}.\]

If $X$ is an object of $\cell(\dtop)$, then a presentation of the
morphism $\varnothing \longrightarrow X$ as a transfinite composition
of pushouts of morphisms of $I^{gl}_+$ is called a \textit{globular
decomposition} of $X$.

\section{Generalized T-homotopy equivalence}
\label{section_defT}

\begin{figure}
\[
\xymatrix{\widehat{0} \ar@{->}[rrrr]^-{U} &&&& \widehat{1} \\
\widehat{0} \ar@{->}[rr]^-{U'} && A \ar@{->}[rr]^-{U''} && \widehat{1}}
\]
\caption{The simplest example of refinement of observation}
\label{ex1}
\end{figure}

We recall here the definition of a T-homotopy equivalence already
given in \cite{1eme} and \cite{2eme}.

\bd A flow $X$ is {\rm loopless} if for any $\alpha\in X^0$, 
the space $\P_{\alpha,\alpha}X$ is empty. \ed

Recall that a flow is a small category without identity morphisms
enriched over a category of topological spaces.  So the preceding
definition is meaningful.

\begin{lem} \label{ordrestate} If a flow $X$ is loopless, then the
  transitive closure of the set \[\{(\alpha,\beta)\in X^0\p X^0\hbox{
    such that }\P_{\alpha,\beta}X\neq\varnothing\}\] induces a partial
  ordering on $X^0$.
\end{lem}

\bpf If $(\alpha,\beta)$ and $(\beta,\alpha)$ with $\alpha\neq \beta$
belong to the transitive closure, then there exists a finite sequence
$(x_1,\dots,x_\ell)$ of elements of $X^0$ with $x_1=\alpha$,
$x_\ell=\alpha$, $\ell>1$ and for any $m$, $\P_{x_m,x_{m+1}}X$ is
non-empty. Consequently, the space $\P_{\alpha,\alpha}X$ is non-empty
because of the existence of the composition law of $X$: contradiction.
\epf

\bd~\footnote{The statement of the definition is slightly different, but equivalent to the statement given in other parts of this series.} 
A {\rm full directed ball} is a flow $\vec{D}$ such that:
\begin{itemize}
\item $\vec{D}$ is loopless (so by Lemma~\ref{ordrestate}, the set $\vec{D}^0$ is equipped with a partial ordering $\leq$)
\item $(\vec{D}^0,\leq)$ is finite bounded
\item for all $(\alpha,\beta)\in \vec{D}^0\p \vec{D}^0$, the topological space 
$\P_{\alpha,\beta}\vec{D}$ is weakly contractible if $\alpha<\beta$, and empty otherwise by definition of $\leq$. 
\end{itemize}
\ed

Let $\vec{D}$ be a full directed ball. Then by Lemma~\ref{ordrestate},
the set $\vec{D}^0$ can be viewed as a finite bounded poset.
Conversely, if $P$ is a finite bounded poset, let us consider the
\textit{flow} $F(P)$ \textit{associated with} $P$: it is of course
defined as the unique flow $F(P)$ such that $F(P)^0=P$ and
$\P_{\alpha,\beta}F(P)=\{u_{\alpha,\beta}\}$ if $\alpha<\beta$ and
$\P_{\alpha,\beta}F(P)=\varnothing$ otherwise. Then $F(P)$ is a full
directed ball and for any full directed ball $\vec{D}$, the two flows
$\vec{D}$ and $F(\vec{D}^0)$ are weakly S-homotopy equivalent.

Let $\vec{E}$ be another full directed ball. Let
$f:\vec{D}\longrightarrow\vec{E}$ be a morphism of flows preserving
the initial and final states. Then $f$ induces a morphism of posets
from $\vec{D}^0$ to $\vec{E}^0$ such that $f(\min
\vec{D}^0)=\min \vec{E}^0$ and $f(\max \vec{D}^0)=\max \vec{E}^0$. Hence 
the following definition:

\bd \label{definitiondeT}
Let $\mathcal{T}$ be the class of morphisms of posets
$f:P_1\longrightarrow P_2$ such that:
\begin{enumerate}
\item The posets $P_1$ and $P_2$ are finite and bounded. 
\item The morphism of posets $f:P_1 \longrightarrow P_2$ is one-to-one; 
in particular, if $x$ and $y$ are two elements of $P_1$ with $x<y$,
then $f(x)<f(y)$.
\item One has $f(\min P_1)=\min P_2$ and  $f(\max P_1)=\max P_2$.
\end{enumerate}
Then a {\rm generalized T-homotopy equivalence} is a morphism of
$\cof(\{Q(F(f)),f\in\mathcal{T}\})$ where $Q$ is the cofibrant
replacement functor of the model category $\dtop$.
\ed

One can choose a \textit{set} of representatives for each isomorphism
class of finite bounded poset. One obtains a \textit{set} of morphisms
$\overline{\mathcal{T}} \subset \mathcal{T}$ such that there is the
equality of classes \[\cof(\{Q(F(f)),f\in\overline{\mathcal{T}}\}) =
\cof(\{Q(F(f)),f\in\mathcal{T}\}).\] By \cite{model3} Proposition~11.5,
the set of morphisms $\{Q(F(f)),f\in\overline{\mathcal{T}}\}$ permits
the small object argument. Thus, the class of morphisms
$\cof(\{Q(F(f)),f\in\mathcal{T}\})$ contains exactly the retracts of
the morphisms of $\cell(\{Q(F(f)),f\in\mathcal{T}\})$ by
\cite{MR99h:55031} Corollary~2.1.15.

The inclusion of posets $\{\widehat{0} < \widehat{1}\} \subset
\{\widehat{0} < A < \widehat{1}\}$ corresponds to the case of
Figure~\ref{ex1}.

\begin{figure}
\begin{center}
\includegraphics[width=9cm]{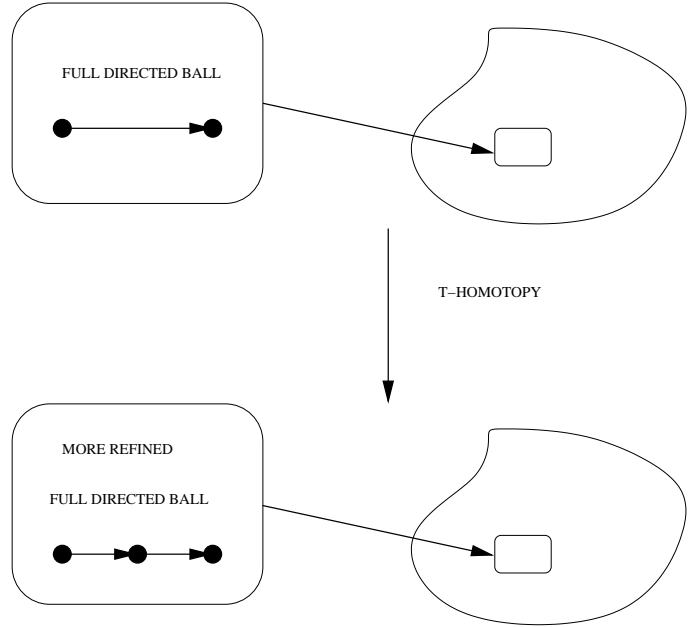}
\end{center}
\caption{Replacement of a full directed ball by a more refined one}
\label{ex2}
\end{figure}

A T-homotopy consists in locally replacing in a flow a full directed
ball by a more refined one (cf. Figure~\ref{ex2}), and in iterating
the process transfinitely.

\section{Principle of the proof of the main theorem}
\label{principlemainproof}

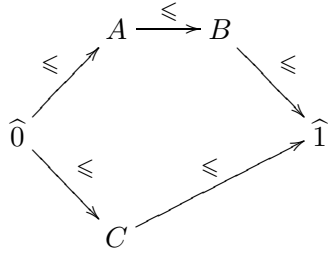
\begin{figure}
\[
\xymatrix{
& A \fr{\leq} & B \ar@{->}[rd]^{\leq} & \\
\widehat{0} \ar@{->}[rd]^{\leq} \ar@{->}[ru]^{\leq} & & & \widehat{1}\\
&  C \ar@{->}[rru]^{\leq} & & }
\] 
\caption{Example of finite bounded poset}
\label{poset0}
\end{figure}

In this section, we collect the main ideas used in the proof of
Theorem~\ref{finpm}. These ideas are illustrated by the case of the
flow $F(P)$ associated with the poset $P$ of Figure~\ref{poset0}.
More precisely, it is going to be explained the reason of the
contractibility of the homotopy branching space
$\hop^-_{\widehat{0}}F(P)$ of the flow $F(P)$ at the initial state
$\widehat{0}$.

First of all, we recall the definition of the branching space functor.
Roughly speaking, the branching space of a flow is the space of germs
of non-constant execution paths beginning in the same way.

\bp \label{universalpm} (\cite{exbranch} Proposition~3.1)
Let $X$ be a flow. There exists a topological space $\P^-X$ unique up
to homeomorphism and a continuous map $h^-:\P X\longrightarrow \P^- X$
satisfying the following universal property:
\begin{enumerate}
\item For any $x$ and $y$ in $\P X$ such that $t(x)=s(y)$, the equality
$h^-(x)=h^-(x*y)$ holds.
\item Let $\phi:\P X\longrightarrow Y$ be a
continuous map such that for any $x$ and $y$ of $\P X$ such that
$t(x)=s(y)$, the equality $\phi(x)=\phi(x*y)$ holds. Then there exists a
unique continuous map $\overline{\phi}:\P^-X\longrightarrow Y$ such that
$\phi=\overline{\phi}\circ h^-$.
\end{enumerate}
Moreover, one has the homeomorphism
\[\P^-X\iso \bigsqcup_{\alpha\in X^0} \P^-_\alpha X\]
where $\P^-_\alpha X:=h^-\lp \bigsqcup_{\beta\in
X^0} \P^-_{\alpha,\beta} X\rp$. The mapping $X\mapsto \P^-X$
yields a functor $\P^-$ from $\dtop$ to $\top$. 
\ep

\bd 
Let $X$ be a flow. The topological space $\P^-X$ is called the {\rm
branching space} of the flow $X$. The functor $\P^-$ is called the 
{\rm branching space functor}. 
\ed

\bth\label{Pleft}  (\cite{exbranch} Theorem~5.5)
The branching space functor \[\P^-:\dtop\longrightarrow \top\] is a left
Quillen functor.
\eth

\bd
The {\rm homotopy branching space} $\hop^- X$ of a flow $X$ is by
definition the topological space $\P^-Q(X)$.  For $\alpha\in X^0$, let
$\hop^-_\alpha X=\P^-_\alpha Q(X)$.
\ed

The first idea would be to replace the calculation of
$\P^-_{\widehat{0}}Q(F(P))$ by the calculation of
$\P^-_{\widehat{0}}F(P)$ because there exists a natural weak
S-homotopy equivalence $Q(F(P)) \longrightarrow F(P)$. However, the
flow $F(P)$ is not cofibrant because its composition law contains
relations, for instance $u_{\widehat{0},A} * u_{A,\widehat{1}} =
u_{\widehat{0},C} * u_{C,\widehat{1}}$.  In any cofibrant replacement
of $F(P)$, a relation like $u_{\widehat{0},A} * u_{A,\widehat{1}} =
u_{\widehat{0},C} * u_{C,\widehat{1}}$ is always replaced by a
S-homotopy between $u_{\widehat{0},A} * u_{A,\widehat{1}}$ and
$u_{\widehat{0},C} * u_{C,\widehat{1}}$. Moreover, it is known after
\cite{exbranch} Theorem~4.1 that the branching space functor does not
necessarily send a weak S-homotopy equivalence of flows to a weak
homotopy equivalence of topological spaces. So this first idea fails,
or at least it cannot work directly.

Let $X=Q(F(P))$ be the cofibrant replacement of $F(P)$. Another idea
that we did not manage to work out can be presented as follows. Every
non-constant execution path $\gamma$ of $\P X$ such that
$s(\gamma)=\widehat{0}$ is in the same equivalence class as an
execution path of $\P_{\widehat{0},\widehat{1}}X$ since the state
$\widehat{1}$ is the only final state of $X$. Therefore, the
topological space $\hop^-_{\widehat{0}}F(P)=\P^-_{\widehat{0}}X$ is a
quotient of the contractible cofibrant space
$\P_{\widehat{0},\widehat{1}}X$. However, the quotient of a
contractible space is not necessarily contractible. For example, 
identifying in the $1$-dimensional disk $\mathbf{D}^1$ the points $-1$
and $+1$ gives the $1$-dimensional sphere $\mathbf{S}^1$.

The principle of the proof exposed in this paper consists in finding a
diagram of topological spaces $\mathcal{P}_{\widehat{0}}^-(X)$
satisfying the following properties:
\begin{enumerate}
\item There is an isomorphism of topological spaces
\[\P^-_{\widehat{0}}X\iso \liminj \mathcal{P}_{\widehat{0}}^-(X).\]  
\item There is a weak homotopy equivalence   of topological spaces 
\[\liminj \mathcal{P}_{\widehat{0}}^-(X)\simeq \holiminj \mathcal{P}_{\widehat{0}}^-(X)\]   
because the diagram of topological spaces
$\mathcal{P}_{\widehat{0}}^-(X)$ is cofibrant for an appropriate model
structure and because for this model structure, the colimit functor is 
a left Quillen functor.
\item Each vertex of the diagram of topological spaces
$\mathcal{P}_{\widehat{0}}^-(X)$ is contractible. Hence, its homotopy
colimit is weakly homotopy equivalent to the classifying space of the
underlying category of $\mathcal{P}_{\widehat{0}}^-(X)$.
\item The underlying category of the diagram
  $\mathcal{P}_{\widehat{0}}^-(X)$ is contractible.
\end{enumerate}

To prove the second assertion, we will build a Reedy structure on the
underlying category of the diagram $\mathcal{P}_{\widehat{0}}^-(X)$.
The main ingredient (but not the only one) of this construction will
be that for every triple $(\alpha,\beta,\gamma)\in X^0\p X^0\p X^0$,
the continuous map $\P_{\alpha,\beta}X\p \P_{\beta,\gamma}X
\longrightarrow \P_{\alpha,\gamma}X$ induced by the composition law of
$X$ is a cofibration of topological spaces since $X$ is cofibrant.

The underlying category of the diagram of topological spaces
$\mathcal{P}_{\widehat{0}}^-(X)$ will be the opposite category 
$\Delta(P\backslash \{\widehat{0}\})^{op} $ of the order complex of the
poset $P\backslash \{\widehat{0}\}$. The latter looks as follows (it is 
the opposite category of the category generated by the inclusions, therefore 
all diagrams are commutative):
\[
\xymatrix{
 & & (A,B,\widehat{1})\ar@{-->}[rd] \ar@{->}[d] \ar@{->}[dl]& \\
(C,\widehat{1})\ar@{-->}[d]\ar@{->}[dr] & \ar@{-->}[dr](A,\widehat{1})\ar@{->}[d] & (B,\widehat{1}) \ar@{-->}[dr]\ar@{->}[dl] & \ar@{-->}[dl](A, B)\ar@{->}[d] \\
(C) & (\widehat{1}) &  (A) & (B)}
\]
The diagram $\mathcal{P}_{\widehat{0}}^-(X)$ is then defined as follows: 
\begin{itemize}
\item $\mathcal{P}_{\widehat{0}}^-(X)(A,B,\widehat{1})=\P_{\widehat{0},A}X\p \P_{A,B}X\p \P_{B,\widehat{1}}X$
\item $\mathcal{P}_{\widehat{0}}^-(X)(A)=\P_{\widehat{0},A}X$
\item $\mathcal{P}_{\widehat{0}}^-(X)(B)=\P_{\widehat{0},B}X$
\item $\mathcal{P}_{\widehat{0}}^-(X)(C)=\P_{\widehat{0},C}X$
\item $\mathcal{P}_{\widehat{0}}^-(X)(\widehat{1})=\P_{\widehat{0},\widehat{1}}X$
\item  $\mathcal{P}_{\widehat{0}}^-(X)(A,B)=\P_{\widehat{0},A}X\p \P_{A,B}X$ 
\item  $\mathcal{P}_{\widehat{0}}^-(X)(B,\widehat{1})=\P_{\widehat{0},B}X\p \P_{B,\widehat{1}}X$
\item  $\mathcal{P}_{\widehat{0}}^-(X)(A,\widehat{1})=\P_{\widehat{0},A}X\p \P_{A,\widehat{1}}X$
\item  $\mathcal{P}_{\widehat{0}}^-(X)(C,\widehat{1})=\P_{\widehat{0},C}X\p \P_{C,\widehat{1}}X$
\item the morphisms $\xymatrix@1{\ar@{-->}[r]&}$ are induced by the projection
\item the morphisms $\xymatrix@1{\ar@{->}[r]&}$ are induced by the composition law.
\end{itemize}

Note that the restriction
$\P_{\widehat{0}}^-(X)\!\restriction_{p_{\widehat{0}}^-(X)}$ of the
diagram of topological spaces $\mathcal{P}_{\widehat{0}}^-(X)$ to the
small category $p_{\widehat{0}}^-(X) \subset
\Delta(P\backslash \{\widehat{0}\})^{op}$
\[
\xymatrix{
(C,\widehat{1})\ar@{-->}[d]\ar@{->}[dr] & \ar@{-->}[dr](A,\widehat{1})\ar@{->}[d] & (B,\widehat{1}) \ar@{-->}[dr]\ar@{->}[dl] & \ar@{-->}[dl](A, B)\ar@{->}[d] \\
(C) & (\widehat{1}) &  (A) & (B)}
\]
has the same colimit, that is $\P_{\widehat{0}}^-(X)$, since the
category $p_{\widehat{0}}^-(X)$ is a final subcategory of
$\Delta(P\backslash \{\widehat{0}\})^{op}$. However, the latter
restriction cannot be Reedy cofibrant because of the associativity of
the composition law. Indeed, the continuous map
\[\P_{\widehat{0},A}X\p \P_{A,\widehat{1}}X\sqcup
\P_{\widehat{0},B}X\p \P_{B,\widehat{1}}X\longrightarrow
\P_{\widehat{0}}^-(X)\!\restriction_{p_{\widehat{0}}^-(X)}(\widehat{1})=\P_{\widehat{0},\widehat{1}}X\]
induced by the composition law of $X$ is not even a monomorphism: if
$(u,v,w)\in \P_{\widehat{0},A}X\p \P_{A,B}X\p \P_{B,\widehat{1}}X$,
then $u*v*w=(u*v)*w\in \P_{\widehat{0},B}X\p \P_{B,\widehat{1}}X$ and
$u*v*w=u*(v*w)\in \P_{\widehat{0},A}X\p \P_{A,\widehat{1}}X$. So the
second assertion of the main argument cannot be true. Moreover, the
classifying space of $p_{\widehat{0}}^-(X)$ is not contractible: it is
homotopy equivalent to the circle $\mathbf{S}^1$. So the fourth
assertion of the main argument cannot be applied either.

On the other hand, the continuous map
\[(\P_{\widehat{0},A}X\p
\P_{A,\widehat{1}}X)\sqcup_{(\P_{\widehat{0},A}X\p \P_{A,B}X\p
  \P_{B,\widehat{1}}X)} (\P_{\widehat{0},B}X\p
\P_{B,\widehat{1}}X)\longrightarrow
\P_{\widehat{0}}^-(X)(\widehat{1})=\P_{\widehat{0},\widehat{1}}X\] is
a cofibration of topological spaces and the classifying space of the
order complex of the poset $P\backslash \{\widehat{0}\}$ is
contractible since the poset $P\backslash \{\widehat{0}\} =
]\widehat{0},\widehat{1}]$ has a unique top element $\widehat{1}$
\cite{MR493916}.

\section{Calculating the  branching space of a loopless flow}
\label{maintool}

\bth  \label{deffunctorP}
Let $X$ be a loopless flow. Let $\alpha\in X^0$. There exists one and
only one functor
\[\mathcal{P}_\alpha^-(X):\Delta(X^0_{>\alpha})^{op}\longrightarrow \top\]
satisfying the following conditions: 
\begin{enumerate}
\item $\mathcal{P}_\alpha^-(X)_{(\alpha_0,\dots,\alpha_p)}:=\P_{\alpha,\alpha_0}X \p \P_{\alpha_0,\alpha_1}X\p\dots\p \P_{\alpha_{p-1},\alpha_p}X$
\item the morphism 
$\de_i:\mathcal{P}_\alpha^-(X)_{(\alpha_0,\dots,\alpha_p)} \longrightarrow
\mathcal{P}_\alpha^-(X)_{(\alpha_0,\dots,\widehat{\alpha_i},\dots,\alpha_p)}$
for $0< i<p$ is induced by the composition law of $X$, more precisely
by the morphism \[\P_{\alpha_{i-1},\alpha_{i}}X \p
\P_{\alpha_{i},\alpha_{i+1}}X \longrightarrow
\P_{\alpha_{i-1},\alpha_{i+1}}X\] 
\item the morphism 
$\de_0:\mathcal{P}_\alpha^-(X)_{(\alpha_0,\dots,\alpha_p)}
\longrightarrow \mathcal{P}_\alpha^-(X)_{(\widehat{\alpha_0},\dots,\alpha_i,\dots,\alpha_p)}$
is induced by the composition law of $X$, more precisely by the
morphism \[\P_{\alpha,\alpha_{0}}X \p
\P_{\alpha_{0},\alpha_{1}}X \longrightarrow
 \P_{\alpha,\alpha_{1}}X\]
\item the morphism 
$\de_p:\mathcal{P}_\alpha^-(X)_{(\alpha_0,\dots,\alpha_p)}
\longrightarrow \mathcal{P}_\alpha^-(X)_{(\alpha_0,\dots,\alpha_{p-1},\widehat{\alpha_p})}$
 is the projection map obtained by removing the component
 $\P_{\alpha_{p-1},\alpha_p}X$.
\end{enumerate}
\eth

\bpf The uniqueness on objects is exactly the first assertion.  The
uniqueness on morphisms comes from the fact that every morphism of
$\Delta(X^0_{>\alpha})^{op}$ is a composite of $\de_i$. We have to
prove the existence.

The diagram  of topological spaces 
\[
\xymatrix{
\mathcal{P}_\alpha^-(X)_{(\alpha_0,\dots,\alpha_p)} \fr{\de_i}\fd{\de_j} & \mathcal{P}_\alpha^-(X)_{(\alpha_0,\dots,\widehat{\alpha_i},\dots,\alpha_p)}\fd{\de_{j-1}}\\
\mathcal{P}_\alpha^-(X)_{(\alpha_0,\dots,\widehat{\alpha_j},\dots,\alpha_p)}\fr{\de_i}& \mathcal{P}_\alpha^-(X)_{(\alpha_0,\dots,\widehat{\alpha_i},\dots,\widehat{\alpha_j},\dots,\alpha_p)}}
\]
is commutative for any $0<i<j<p$ and any $p\geq 2$. Indeed, if $i<j-1$, 
then one has
\[\boxed{\de_i\de_j(\gamma_0,\dots,\gamma_{p})=\de_{j-1}\de_i(\gamma_0,\dots,\gamma_{p})=
(\gamma_0,\dots,\gamma_{i}\gamma_{i+1},\dots,\gamma_{j}\gamma_{j+1},\dots,\gamma_{p})}\]
and if $i=j-1$, then one has 
\[\boxed{\de_i\de_j(\gamma_0,\dots,\gamma_{p})=\de_{j-1}\de_i(\gamma_0,\dots,\gamma_{p})=
(\gamma_0,\dots,\gamma_{j-1}\gamma_{j}\gamma_{j+1},\dots,\gamma_{p})}\]
because of the associativity of the composition law of $X$ (it is
the only place in this proof where this axiom is required).

The diagram  of topological spaces 
\[
\xymatrix{
\mathcal{P}_\alpha^-(X)_{(\alpha_0,\dots,\alpha_p)} \fr{\de_i}\fd{\de_p} & \mathcal{P}_\alpha^-(X)_{(\alpha_0,\dots,\widehat{\alpha_i},\dots,\alpha_p)}\fd{\de_{p-1}}\\
\mathcal{P}_\alpha^-(X)_{(\alpha_0,\dots,\alpha_{p-1})}\fr{\de_i}& \mathcal{P}_\alpha^-(X)_{(\alpha_0,\dots,\widehat{\alpha_i},\dots,\alpha_{p-1})}}
\]
is commutative for any $0<i<p-1$ and any $p> 2$. Indeed, one has
\[\boxed{\de_i\de_p(\gamma_0,\dots,\gamma_{p})=\de_{p-1}\de_i(\gamma_0,\dots,\gamma_{p})=(\gamma_0,\dots,\gamma_{i}\gamma_{i+1},\dots,\gamma_{p-1})}.\]

At last, the diagram  of topological spaces 
\[
\xymatrix{
\mathcal{P}_\alpha^-(X)_{(\alpha_0,\dots,\alpha_p)} \fr{\de_{p-1}}\fd{\de_{p}} & \mathcal{P}_\alpha^-(X)_{(\alpha_0,\dots,\widehat{\alpha_{p-1}},\alpha_p)}\fd{\de_{p-1}}\\
\mathcal{P}_\alpha^-(X)_{(\alpha_0,\dots,\alpha_{p-1})}\fr{\de_{p-1}}& \mathcal{P}_\alpha^-(X)_{(\alpha_0,\dots,\alpha_{p-2})}}
\]
is commutative for any $p\geq 2$.  Indeed, one has
\[\boxed{\de_{p-1}\de_p(\gamma_0,\dots,\gamma_{p})=(\gamma_0,\dots,\gamma_{p-2})}\]
and 
\[\boxed{\de_{p-1}\de_{p-1}(\gamma_0,\dots,\gamma_{p})=\de_{p-1}(\gamma_0,\dots,\gamma_{p-2},\gamma_{p-1}\gamma_{p})=(\gamma_0,\dots,\gamma_{p-2})}.\]

In other terms, the $\de_i$ maps satisfy the simplicial
identities. Hence the result. \epf

The following theorem is used in the proofs of Theorem~\ref{calculPm}
and Theorem~\ref{simplification}.

\bth (\cite{MR1712872} Theorem~1 p.~213) \label{souslim}
Let $L:J'\longrightarrow J$ be a final functor between small
categories, i.e.  such that for any $k\in J$, the comma category
$(k\!\downarrow\!L)$ is non-empty and connected. Let
$F:J\longrightarrow \C$ be a functor from $J$ to a cocomplete category
$\C$. Then $L$ induces a canonical morphism $\liminj (F\circ
L)\longrightarrow \liminj F$ which is an isomorphism.
\eth

\bth \label{calculPm} 
Let $X$ be a loopless flow. Then there exists an isomorphism of
topological spaces $\P_\alpha^- X\iso \liminj
\mathcal{P}_\alpha^-(X)$ for any $\alpha\in X^0$.
\eth

\bpf
Let $p_\alpha^-(X)$ be the full subcategory of
$\Delta(X^0_{>\alpha})^{op}$ generated by the arrows
$\de_0:(\alpha_0,\alpha_1)\longrightarrow (\alpha_1)$ and
$\de_1:(\alpha_0,\alpha_1)\longrightarrow (\alpha_0)$.

Let $k=(k_0,\dots,k_q)$ be an object of
$\Delta(X^0_{>\alpha})^{op}$. Then $k\to (k_0)$ is an object of the comma
category $(k\!\downarrow\!p_\alpha^-(X))$. So the latter category is
not empty. Let $k\to (x_0)$ and $k\to (y_0)$ be two distinct elements
of $(k\!\downarrow\!p_\alpha^-(X))$. The pair $\{x_0,y_0\}$ is
therefore a subset of $\{k_0,\dots,k_q\}$. So either $x_0<y_0$ or
$y_0<x_0$. Without loss of generality, one can suppose that
$x_0<y_0$. Then one has the commutative diagram
\[
\xymatrix{
k\ar@{=}[r]\fd{} & k\ar@{=}[r]\fd{} & k\fd{} \\
(x_0) & (x_0,y_0)\fl{\de_1}\fr{\de_0} & (y_0).}
\] 
Therefore, the objects $k\to (x_0)$ and $k\to (y_0)$ are in the same
connected component of $(k\!\downarrow\!p_\alpha^-(X))$.  Let $k\to
(x_0)$ and $k\to (y_0,y_1)$ be two distinct elements of
$(k\!\downarrow\!p_\alpha^-(X))$. Then $k\to (x_0)$ is in the same
connected component as $k\to (y_0)$ by the previous calculation.
Moreover, one has the commutative diagram
\[
\xymatrix{
k\ar@{=}[r]\fd{} & k\fd{} \\
(y_0) & (y_0,y_1).\fl{\de_1}}
\] 
Thus, the objects $k\to (x_0)$ and $k\to (y_0,y_1)$ are in the same
connected component of $(k\!\downarrow\!p_\alpha^-(X))$. So the comma
category $(k\!\downarrow\!p_\alpha^-(X))$ is connected and non-empty.
Thus for any functor $F:\Delta(X^0_{>\alpha})^{op}\longrightarrow
\top$, the inclusion functor $i:p_\alpha^-(X)\longrightarrow
\Delta(X^0_{>\alpha})^{op}$ induces an isomorphism of topological
spaces $\liminj (F\circ i)\longrightarrow\liminj F$ by
Theorem~\ref{souslim}.

Let $\widehat{p}_\alpha^-(X)$ be the full subcategory of
$p_\alpha^-(X)$ consisting of the objects $(\alpha_0)$. The category
$\widehat{p}_\alpha^-(X)$ is discrete because it does not contain any
non identity morphism. Let $j:\widehat{p}_\alpha^-(X)\longrightarrow
p_\alpha^-(X)$ be the canonical inclusion functor. It induces a
canonical continuous map $\liminj (F\circ j)\longrightarrow\liminj
(F\circ i)$ for any functor
$F:\Delta(X^0_{>\alpha})^{op}\longrightarrow \top$.

For $F=\mathcal{P}_\alpha^-(X)$, one obtains the diagram of topological spaces 
\[\liminj (\mathcal{P}_\alpha^-(X)\circ j)\longrightarrow\liminj (\mathcal{P}_\alpha^-(X)\circ i)\iso \liminj \mathcal{P}_\alpha^-(X).\]
It is clear that $\liminj (\mathcal{P}_\alpha^-(X)\circ
j)\iso\bigsqcup_{\alpha_0}\P_{\alpha,\alpha_0}X$. Let $g:\liminj
(\mathcal{P}_\alpha^-(X)\circ j)\longrightarrow Z$ be a continuous map such
that $g(x*y)=g(x)$ for any $x$ and any $y$ such that $t(x)=s(y)$. So
there exists a commutative diagram
\[
\xymatrix{
\P_{s(x),t(x)}X\p \P_{t(x),t(y)}X \fd{\de_1}\fr{\de_0}\ar@{->}[rd] & \P_{s(x),t(y)}X \fd{g}\\
\P_{s(x),t(x)}X \fr{g}& Z}
\]
for any $x$ and $y$ as above. Therefore, the topological space
$\liminj (\mathcal{P}_\alpha^-(X)\circ i)$ satisfies the same
universal property as the topological space $\P_\alpha^-X$ (cf.
Proposition~\ref{universalpm}).  \epf

\section{Reedy structure and homotopy colimit}
\label{reedythree}

\begin{lem}
  Let $X$ be a loopless flow such that $(X^0,\leq)$ is locally finite.
  If $(\alpha,\beta)$ is a $1$-simplex of $\Delta(X^0)$ and if
  $(\alpha_0,\dots,\alpha_p)$ is a $p$-simplex of $\Delta(X^0)$ with
  $\alpha_0=\alpha$ and $\alpha_p=\beta$, then $p$ is at most the
  cardinal $\card(]\alpha,\beta])$ of $]\alpha,\beta]$.
\end{lem}

\bpf If $(\alpha_0,\dots,\alpha_p)$ is a $p$-simplex of $\Delta(X^0)$,
then one has $\alpha_0<\dots<\alpha_p$ by definition of the order
complex. So one has the inclusion $\{\alpha_1,\dots,\alpha_p\}\subset
]\alpha,\beta]$, and therefore $p\leq \card(]\alpha,\beta])$. \epf

The following choice of notation is therefore meaningful.

\begin{nota} \label{defell} 
Let $X$ be a loopless flow such that $(X^0,\leq)$ is locally finite. 
Let $(\alpha,\beta)$ be a $1$-simplex of $\Delta(X^0)$. We denote by
$\ell(\alpha,\beta)$ the maximum of the set of integers
\[\left\{p\geq 1, \exists (\alpha_0,\dots,\alpha_p)\hbox{ $p$-simplex of }\Delta(X^0)\hbox{ s.t. }
(\alpha_0,\alpha_p)= (\alpha,\beta)\right\}\] One always has $1\leq
\ell(\alpha,\beta)\leq \card(]\alpha,\beta])$.
\end{nota}

\begin{lem}\label{triinv} 
Let $X$ be a loopless flow such that $(X^0,\leq)$ is locally finite.
Let $(\alpha,\beta,\gamma)$ be a $2$-simplex of $\Delta(X^0)$. Then
one has
\[\ell(\alpha,\beta)+\ell(\beta,\gamma)\leq \ell(\alpha,\gamma).\] 
\end{lem}

\bpf 
Let $\alpha = \alpha_0 < \dots < \alpha_{\ell(\alpha,\beta)} =
\beta$. Let $\beta = \beta_0 < \dots < \beta_{\ell(\beta,\gamma)} = 
\gamma$. Then 
\[(\alpha_0,\dots,\alpha_{\ell(\alpha,\beta)},\beta_1,\dots,\beta_{\ell(\beta,\gamma)}) \]
is a simplex of $\Delta(X^0)$ with $\alpha = \alpha_0$ and $\beta_{\ell(\beta,\gamma)} = 
\gamma$. So $\ell(\alpha,\beta)+\ell(\beta,\gamma)\leq \ell(\alpha,\gamma)$. 
\epf

\bp \label{ReedyOne}
Let $X$ be a loopless flow such that $(X^0,\leq)$ is locally
finite. Let $\alpha\in X^0$. Let $\Delta(X^0_{>\alpha})^{op}_+$ be the
subcategory of $\Delta(X^0_{>\alpha})^{op}$ generated by the
\[\de_i:(\alpha_0,\dots,\alpha_p)\longrightarrow
(\alpha_0,\dots,\widehat{\alpha_i},\dots,\alpha_p)\] 
for any $p\geq 1$ and $0\leq i<p$. Let $\Delta(X^0_{>\alpha})^{op}_-$ be the
subcategory of $\Delta(X^0_{>\alpha})^{op}$ generated by the 
\[\de_p:(\alpha_0,\dots,\alpha_p)\longrightarrow
(\alpha_0,\dots,\alpha_{p-1})\] 
for any $p\geq 1$. If $(\alpha_0,\dots,\alpha_p)$ 
is an object of $\Delta(X^0_{>\alpha})^{op}$, let: 
\[d(\alpha_0,\dots,\alpha_p) =\ell(\alpha,\alpha_0)^2+\ell(\alpha_0,\alpha_1)^2+\dots+\ell(\alpha_{p-1},\alpha_p)^2.\]
Then the triple
$(\Delta(X^0_{>\alpha})^{op},\Delta(X^0_{>\alpha})^{op}_+,\Delta(X^0_{>\alpha})^{op}_-)$
together with the degree function $d$ is a Reedy category. 
\ep

Note that the subcategory $\Delta(X^0_{>\alpha})^{op}_+$ is precisely
generated by the morphisms sent by the functor
$\mathcal{P}_\alpha^-(X)$ (Theorem~\ref{deffunctorP}) to continuous
maps induced by the composition law of the flow $X$. And note that the
subcategory $\Delta(X^0_{>\alpha})^{op}_-$ is precisely generated by
the morphisms sent by the functor $\mathcal{P}_\alpha^-(X)$
(Theorem~\ref{deffunctorP}) to continuous maps induced by the
projection obtained by removing the last component on the right.

\bpf 
Let $\de_i:(\alpha_0,\dots,\alpha_p)\longrightarrow
(\alpha_0,\dots,\widehat{\alpha_i},\dots,\alpha_p)$ be a morphism of
$\Delta(X^0_{>\alpha})^{op}_+$ with $p\geq 1$ and $0\leq i<p$ .  Then
(with the convention $\alpha_{-1}=\alpha$)
\beas
&&d(\alpha_0,\dots,\alpha_p) =\ell(\alpha,\alpha_0)^2+\dots+\ell(\alpha_{p-1},\alpha_p)^2\\
&& d(\alpha_0,\dots,\widehat{\alpha_i},\dots,\alpha_p) =\ell(\alpha,\alpha_0)^2+\dots+\ell(\alpha_{i-1},\alpha_{i+1})^2+\dots+\ell(\alpha_{p-1},\alpha_p)^2.
\eeas
So one obtains 
\[ d(\alpha_0,\dots,\alpha_p)-d(\alpha_0,\dots,\widehat{\alpha_i},\dots,\alpha_p)
=\ell(\alpha_{i-1},\alpha_{i})^2+\ell(\alpha_{i},\alpha_{i+1})^2-\ell(\alpha_{i-1},\alpha_{i+1})^2.\]
By Lemma~\ref{triinv}, one has 
\[(\ell(\alpha_{i-1},\alpha_{i})+\ell(\alpha_{i},\alpha_{i+1}))^2\leq 
\ell(\alpha_{i-1},\alpha_{i+1})^2\]
and one has 
\[\ell(\alpha_{i-1},\alpha_{i})^2+\ell(\alpha_{i},\alpha_{i+1})^2<(\ell(\alpha_{i-1},\alpha_{i})+\ell(\alpha_{i},\alpha_{i+1}))^2
\]
since $2\ell(\alpha_{i-1},\alpha_{i})\ell(\alpha_{i},\alpha_{i+1})\geq
2$. So every morphism of $\Delta(X^0_{>\alpha})^{op}_+$ raises the degree.

Let $\de_p:(\alpha_0,\dots,\alpha_p)\longrightarrow
(\alpha_0,\dots,\alpha_{p-1})$ with $p\geq 1$ be a morphism of
$\Delta(X^0_{>\alpha})^{op}_-$. Then one has
\beas
&&d(\alpha_0,\dots,\alpha_p) =\ell(\alpha,\alpha_0)^2+\dots+\ell(\alpha_{p-2},\alpha_{p-1})^2+\ell(\alpha_{p-1},\alpha_p)^2\\
&& d(\alpha_0,\dots,\alpha_{p-1}) =\ell(\alpha,\alpha_0)^2+\dots+\ell(\alpha_{p-2},\alpha_{p-1})^2.
\eeas
So
$d(\alpha_0,\dots,\alpha_p)-d(\alpha_0,\dots,\alpha_{p-1}) = \ell(\alpha_{p-1},\alpha_p)^2>0$.
Thus, every morphism of $\Delta(X^0_{>\alpha})^{op}_-$ lowers the degree.

Let $f:(\alpha_0,\dots,\alpha_p)\longrightarrow
(\beta_0,\dots,\beta_q)$ be a morphism of
$\Delta(X^0_{>\alpha})^{op}$.  Then one has $p\geq q$ and
$(\beta_0,\dots,\beta_q)=(\alpha_{\sigma(0)},\dots,\alpha_{\sigma(q)})$
where $\sigma:\{0,\dots,q\}\longrightarrow \{0,\dots,p\}$ is a
strictly increasing set map. Then $f$ can be written as a composite
\[\xymatrix{(\alpha_0,\alpha_1,\alpha_2,\dots,\alpha_p)\ar@{->}[rr]^-{}&& (\alpha_0,\alpha_1,\alpha_2,\dots,\alpha_{\sigma(q)})\ar@{->}[rr]^-{\de_{j_1}\de_{j_2}\dots\de_{j_r}}&&(\beta_1,\dots,\beta_q)}\]
where $0\leq j_1<j_2<\dots<j_r<\sigma(q)$ and
$\{j_1,j_2,\dots,j_r\}\cup
\{\sigma(1),\sigma(2),\dots,\sigma(q)\}=\{0,1,2,\dots,\sigma(q)\}$ and
this is the unique way of decomposing $f$ as a morphism of
$\Delta(X^0_{>\alpha})^{op}_-$ followed by a morphism of
$\Delta(X^0_{>\alpha})^{op}_+$.  \epf

\bth\label{equi} 
Let $X$ be a loopless flow such that $(X^0,\leq)$ is locally
finite. Let $\alpha\in X^0$. Then the colimit functor 
\[\liminj:\top^{\Delta(X^0_{>\alpha})^{op}}\longrightarrow \top\] 
is a left Quillen functor if the category of diagrams
$\top^{\Delta(X^0_{>\alpha})^{op}}$ is equipped with the Reedy model
structure.
\eth

\bpf The Reedy structure on $\Delta(X^0_{>\alpha})^{op}$ provides a model structure 
on the category $\top^{\Delta(X^0_{>\alpha})^{op}}$ of diagrams of
topological spaces over $\Delta(X^0_{>\alpha})^{op}$ such that a morphism
of diagrams $f:D\longrightarrow E$ is
\begin{enumerate}
\item a weak equivalence if and only if for all objects
  ${\underline{\alpha}}$ of $\Delta(X^0_{>\alpha})^{op}$, the morphism
  $D_{\underline{\alpha}}\longrightarrow E_{\underline{\alpha}}$ is a
  weak homotopy equivalence of $\top$, that is an objectwise weak
  homotopy equivalence
\item a cofibration if and only if for all objects ${\underline{\alpha}}$ of 
$\Delta(X^0_{>\alpha})^{op}$, the morphism
$D_{\underline{\alpha}}\sqcup_{L_{\underline{\alpha}} D}
L_{\underline{\alpha}}E \longrightarrow E_{\underline{\alpha}}$ is a
cofibration of $\top$
\item a fibration if and only if for all objects ${\underline{\alpha}}$
  of $\Delta(X^0_{>\alpha})^{op}$, the morphism
  $D_{\underline{\alpha}}\longrightarrow
  E_{\underline{\alpha}}\p_{M_{\underline{\alpha}} E}
  M_{\underline{\alpha}}D$ is a fibration of $\top$.
\end{enumerate}
Consider the categorical adjunction
$\liminj:\top^{\Delta(X^0_{>\alpha})^{op}}\leftrightarrows \top:\diag$
where $\diag$ is the diagonal functor. Let $p:E\longrightarrow B$ be a
continuous map.  If ${\underline{\alpha}}=(\alpha_0)$, then the
matching category of ${\underline{\alpha}}$ is the empty category. So
$\diag(B)_{\underline{\alpha}}\p_{M_{\underline{\alpha}} \diag(B)}
M_{\underline{\alpha}}\diag(E)\iso \diag(B)_{\underline{\alpha}}\iso
B$ since the spaces $M_{\underline{\alpha}} \diag(B)$ and
$M_{\underline{\alpha}}\diag(E)$ are both equal to singletons. If
${\underline{\alpha}}=(\alpha_0,\dots,\alpha_p)$ with $p\geq 1$, then
the matching category of ${\underline{\alpha}}$ looks like the
following tower:
\[(\alpha_0,\dots,\alpha_{p-1})\longrightarrow (\alpha_0,\dots,\alpha_{p-2})
\longrightarrow\dots \longrightarrow(\alpha_0).\]
Therefore in that case, one has the isomorphisms \[M_{\underline{\alpha}}
\diag(B)\iso \diag(B)_{(\alpha_0,\dots,\alpha_{p-1})}
\iso B\] and \[M_{\underline{\alpha}} \diag(E)\iso \diag(E)_{(\alpha_0,\dots,\alpha_{p-1})}
\iso E.\] Hence $\diag(B)_{\underline{\alpha}}\p_{M_{\underline{\alpha}} \diag(B)}
M_{\underline{\alpha}}\diag(E)\iso E$. Thus, the continuous map
\[E\iso \diag(E)_{\underline{\alpha}}\longrightarrow
\diag(B)_{\underline{\alpha}}\p_{M_{\underline{\alpha}} \diag(B)}
M_{\underline{\alpha}}\diag(E)\] is either the identity of $E$, or
$p$. So if $p$ is a fibration (resp. a trivial fibration), then
$\diag(p):\diag(E)\longrightarrow \diag(B)$ is a Reedy fibration
(resp. a Reedy trivial fibration). One then deduces that $\diag$ is a
right Quillen functor and that the colimit functor is a left Quillen
functor. \epf

\section{Homotopy branching space  of a full directed ball}
\label{diagcofibrant}

\begin{nota} 
  Let $X$   be a loopless flow. Let
$\underline{\alpha}=(\alpha_0,\dots,\alpha_p)$ and
$\underline{\beta}=(\beta_0,\dots,\beta_q)$ be two simplices of the
order complex $\Delta(X^0)$ of the poset $X^0$. Then we denote by:
\begin{enumerate}
\item $\underline{\beta}\supseteqq \underline{\alpha}$ the following 
situation: 
$\alpha_0=\beta_0$, $\alpha_p=\beta_q$ and
$\{\alpha_0,\dots,\alpha_p\}\subset \{\beta_0,\dots,\beta_q\}$ 
\item  $\underline{\beta}\supsetneqq\underline{\alpha}$ the following situation:
$\alpha_0=\beta_0$, $\alpha_p=\beta_q$ and
$\{\alpha_0,\dots,\alpha_p\}\subset \{\beta_0,\dots,\beta_q\}$ and 
$\{\alpha_0,\dots,\alpha_p\}\neq \{\beta_0,\dots,\beta_q\}$.
\end{enumerate}
\end{nota}

\begin{nota} Let $X$ be a loopless flow. 
Let $\underline{\alpha}=(\alpha_0,\dots,\alpha_p)$ be a simplex of the
order complex $\Delta(X^0)$ of the poset $X^0$. Let
$\alpha<\alpha_0$. Then the notation $\alpha.\underline{\alpha}$
represents the simplex $(\alpha,\alpha_0,\dots,\alpha_p)$ of
$\Delta(X^0)$.
\end{nota}

\bth\label{simplification} 
Let $X$ be a loopless flow such that $(X^0,\leq)$ is locally
finite. Let $\alpha\in X^0$. Then:
\begin{enumerate}
\item For any object ${\underline{\alpha}}$ of
$\Delta(X^0_{>\alpha})^{op}$, one has the isomorphism of topological
spaces
\[
{L_{\underline{\alpha}}\mathcal{P}_\alpha^-(X)\iso\liminj_{{\alpha.\underline{\beta}}\supsetneqq{\alpha.\underline{\alpha}}}\mathcal{P}_\alpha^-(X)_{\underline{\beta}}}. 
\]
\item For any object ${\underline{\alpha}}=(\alpha_0,\dots,\alpha_p)$, 
the canonical continuous map
\[i_{\alpha.\underline{\alpha}}:L_{\underline{\alpha}}\mathcal{P}_\alpha^-(X)\longrightarrow
\mathcal{P}_\alpha^-(X)_{\underline{\alpha}}\] is equal to the pushout
product (cf. Notation~\ref{notapush}) of the canonical continuous maps
\[{i_{(\alpha,\alpha_0)}\square i_{(\alpha_0,\alpha_1)} \square \dots \square i_{(\alpha_{p-1},\alpha_p)}}.\]  
\end{enumerate}
\eth

\bpf Let ${\underline{\alpha}}$ be a fixed object of
$\Delta(X^0_{>\alpha})^{op}$. Since the subcategory
$\Delta(X^0_{>\alpha})^{op}_+$ contains only commutative diagrams, the
latching category $\de(\Delta(X^0_{>\alpha})^{op}_+\!\downarrow\!
{\underline{\alpha}})$ is the full subcategory of
$\Delta(X^0_{>\alpha})^{op}_+$ consisting of the simplices
${\underline{\beta}}$ such that
${\alpha.\underline{\beta}}\supsetneqq{\alpha.\underline{\alpha}}$.
Hence the first assertion.

Let $\boxed{\alpha_{-1}:=\alpha}$. One has
${\alpha.\underline{\beta}}\supsetneqq{\alpha.\underline{\alpha}}$ if and only
if the simplex ${\alpha.\underline{\beta}}$ can be written as an expression of the form 
\[\alpha_{-1}.\underline{\delta_0}.\underline{\delta_1}\dots \underline{\delta_p}\]
with
$\alpha_{i-1}.\underline{\delta_i} \supseteqq (\alpha_{i-1},\alpha_{i})$ for all 
$0\leq i\leq p$ and such that at least for one $i$, one has
$\alpha_{i-1}.\underline{\delta_i}  \supsetneqq (\alpha_{i-1},\alpha_{i})$.

Let $\mathcal{E}$ be the set of subsets $S$ of $\{0,\dots,p\}$ such
that $S\neq\{0,\dots,p\}$. For such a $S$, let $I(S)$ be the full
subcategory of $\Delta(X^0_{>\alpha})^{op}_+$ consisting of the
objects ${\underline{\beta}}$ such that
\begin{itemize}
\item ${\alpha.\underline{\beta}} \supsetneqq {\alpha.\underline{\alpha}}$ 
\item for each $i\notin S$, one has
$\alpha_{i-1}.\underline{\delta_i} \supsetneqq (\alpha_{i-1},\alpha_{i})$, and therefore $\underline{\delta_i} \neq (\alpha_i)$ 
\item for each $i\in S$, one has
$\alpha_{i-1}.\underline{\delta_i} \supseteqq (\alpha_{i-1},\alpha_{i})$. 
\end{itemize}

The full subcategory $\bigcup_{S\in \mathcal{E}} I(S)$ is exactly the
subcategory of $\Delta(X^0_{>\alpha})^{op}_+$ consisting of the objects
${\underline{\beta}}$ such that
${\alpha.\underline{\beta}} \supsetneqq {\alpha.\underline{\alpha}}$, that is to
say the subcategory calculating
$L_{\underline{\alpha}}\mathcal{P}_\alpha^-(X)$. In other terms, one
obtains the isomorphism
\be \label{equ1}
\liminj_{\bigcup_{S\in \mathcal{E}} I(S)}\mathcal{P}_\alpha^-(X)\iso L_{\underline{\alpha}}\mathcal{P}_\alpha^-(X).\ee

The full subcategory $I(S)$ of $\Delta(X^0_{>\alpha})^{op}_+$ has a final
subcategory $\overline{I(S)}$ consisting of the objects $\underline{\beta}$
such that
\begin{itemize}
\item ${\alpha.\underline{\beta}} \supsetneqq {\alpha.\underline{\alpha}}$ 
\item for each $i\notin S$, one has
$\alpha_{i-1}.\underline{\delta_i} \supsetneqq (\alpha_{i-1},\alpha_{i})$, and therefore $\underline{\delta_i} \neq (\alpha_i)$ 
\item for each $i\in S$, one has
$\alpha_{i-1}.\underline{\delta_i} =  (\alpha_{i-1},\alpha_{i})$ and therefore $\underline{\delta_i} = (\alpha_i)$.
\end{itemize}
The subcategory $\overline{I(S)}$ is final in $I(S)$ because for
any object $\underline{\beta}$ of $I(S)$, there exists a unique
$\underline{\gamma}$ of $\overline{I(S)}$ and a unique arrow
$\underline{\beta}\longrightarrow
\underline{\gamma}$. Therefore, by Theorem~\ref{souslim}, there is an isomorphism
\be \label{equ2}
\liminj_{I(S)}\mathcal{P}_\alpha^-(X)\iso \liminj_{\overline{I(S)}}\mathcal{P}_\alpha^-(X)
\ee
since the comma category
$(\underline{\beta}\!\downarrow\!\overline{I(S)})$ is the one-object 
category.

For any object $\underline{\beta}$ of $\overline{I(S)}$, one gets  
\begin{align*}
& \mathcal{P}_\alpha^-(X)_{\underline{\beta}} & \\
& =
\prod_{i=0}^{i=p} \mathcal{P}_{\alpha_{i-1}}^-(X)_{\underline{\delta_i}} & \hbox{ by definition of $\mathcal{P}^-$} \\
& \iso \lp \prod_{i \in S}\mathcal{P}_{\alpha_{i-1}}^-(X)_{(\alpha_{i})}\rp 
\p \lp \prod_{i \notin S} \mathcal{P}_{\alpha_{i-1}}^-(X)_{\underline{\delta_i}}\rp 
& \hbox{ by definition of $S$.}
\end{align*}
Thus, since the category $\top$ of compactly generated topological
spaces is cartesian closed, one obtains 
{\small \begin{align*}
&
\liminj_{\overline{I(S)}} \mathcal{P}_\alpha^-(X) & \\
& \iso \liminj_{\overline{I(S)}} \lp\lp \prod_{i \in S}\mathcal{P}_{\alpha_{i-1}}^-(X)_{(\alpha_{i})}\rp 
\p \lp \prod_{i \notin S} \mathcal{P}_{\alpha_{i-1}}^-(X)_{\underline{\delta_i}}\rp\rp & \\
& \iso \lp \prod_{i \in S}\mathcal{P}_{\alpha_{i-1}}^-(X)_{(\alpha_{i})}\rp 
\p \liminj_{\begin{array}{c} i \notin S \\ \alpha_{i-1}.\underline{\delta_i}  \supsetneqq (\alpha_{i-1},\alpha_{i}) \end{array}} \lp \prod_{i \notin S} \mathcal{P}_{\alpha_{i-1}}^-(X)_{\underline{\delta_i}}\rp 
& \\
& \iso \lp \prod_{i \in S}\mathcal{P}_{\alpha_{i-1}}^-(X)_{(\alpha_{i})}\rp 
\p \lp \prod_{i \notin S} \liminj_{\alpha_{i-1}.\underline{\delta_i}  \supsetneqq (\alpha_{i-1},\alpha_{i})}\mathcal{P}_{\alpha_{i-1}}^-(X)_{\underline{\delta_i}}\rp & 
\hbox{ by Lemma~\ref{colimproduit}}.
\end{align*}}
Therefore, one obtains the isomorphism of topological spaces 
\be \label{equ3}
\liminj_{\overline{I(S)}} \mathcal{P}_\alpha^-(X)
\iso \lp \prod_{i \in
  S}\mathcal{P}_{\alpha_{i-1}}^-(X)_{(\alpha_{i})}\rp \p \lp \prod_{i
  \notin S} L_{(\alpha_{i})} \mathcal{P}_{\alpha_{i-1}}^-(X) \rp 
\ee
thanks to the first assertion of the theorem.

If $S$ and $T$ are two elements of $\mathcal{E}$ such that
$S\subset T$, then there exists a canonical morphism of diagrams
$I(S)\longrightarrow I(T)$ inducing a canonical morphism of
topological spaces
\[\liminj_{\underline{\beta}\in
I(S)}\mathcal{P}^-(X)_{\underline{\beta}}
\longrightarrow  \liminj_{\underline{\beta}\in
I(T)}\mathcal{P}^-(X)_{\underline{\beta}}.\] 
Therefore, by Equation~(\ref{equ2}) and Equation~(\ref{equ3}), the double colimit 
\[\liminj_{S\in \mathcal{E}}\lp \liminj_{I(S)} \mathcal{P}_\alpha^-(X)\rp\]
calculates the source of the morphism ${i_{(\alpha,\alpha_0)}\square
i_{(\alpha_0,\alpha_1)} \square \dots \square
i_{(\alpha_{p-1},\alpha_p)}}$ by Theorem~\ref{calculpushout}.  It
then suffices to prove the isomorphism
\[\liminj_{S\in \mathcal{E}}\lp \liminj_{I(S)} \mathcal{P}_\alpha^-(X)\rp \iso  
\liminj_{{\alpha.\underline{\beta}}\supsetneqq{\alpha.\underline{\alpha}}}\mathcal{P}_\alpha^-(X)_{\underline{\beta}}\] 
to complete the proof. For that purpose, it suffices to construct two
canonical morphisms
\[\liminj_{S\in \mathcal{E}}\lp \liminj_{I(S)} \mathcal{P}_\alpha^-(X)\rp 
\longrightarrow \liminj_{{\alpha.\underline{\beta}} \supsetneqq {\alpha.\underline{\alpha}}}\mathcal{P}_\alpha^-(X)_{\underline{\beta}}\] 
and 
\[\liminj_{{\alpha.\underline{\beta}} \supsetneqq {\alpha.\underline{\alpha}}}\mathcal{P}_\alpha^-(X)_{\underline{\beta}} \longrightarrow \liminj_{S\in \mathcal{E}}\lp \liminj_{I(S)} \mathcal{P}_\alpha^-(X)\rp.\] 
The first morphism comes from the isomorphism of
Equation~(\ref{equ1}). As for the second morphism, let us consider a
diagram of flows of the form:
\[
\xymatrix{
\mathcal{P}_\alpha^-(X)_{\underline{\beta}} \fd{} \fr{} & \liminj_{S\in \mathcal{E}}\lp \liminj_{I(S)} \mathcal{P}_\alpha^-(X)\rp\\
\mathcal{P}_\alpha^-(X)_{\underline{\gamma}} \ar@{->}[ru] &}
\]
One has to prove that it is commutative. Since one has $\bigcup_{S\in
\mathcal{E}} I(S)=\Delta(X^0_{>\alpha})^{op}_+$, there exists $S\in \mathcal{E}$
such that $\underline{\gamma}$ is an object of $I(S)$. So
$\underline{\beta}$ is an object of $I(S)$ as well and there exists a
commutative diagram
\[
\xymatrix{
\mathcal{P}_\alpha^-(X)_{\underline{\beta}} \fd{} \fr{} &  \liminj_{I(S)} \mathcal{P}_\alpha^-(X)\\
\mathcal{P}_\alpha^-(X)_{\underline{\gamma}} \ar@{->}[ru] &}
\]
since the subcategory $\Delta(X^0_{>\alpha})^{op}_+$ is
commutative. Hence the result.
\epf

\bth \label{cascofibrant} 
Let $X$ be a loopless object of $\hda$ such that $(X^0,\leq)$ is
locally finite. Let $\alpha\in X^0$. Then the diagram of topological
spaces $\mathcal{P}_\alpha^-(X)$ is Reedy cofibrant. In other terms,
for any object ${\underline{\alpha}}$ of $\Delta(X^0_{>\alpha})^{op}$,
the topological space $\mathcal{P}_\alpha^-(X)_{\underline{\alpha}}$
is cofibrant and the morphism
$L_{\underline{\alpha}}\mathcal{P}_\alpha^-(X)\longrightarrow
\mathcal{P}_\alpha^-(X)_{\underline{\alpha}}$ 
is a cofibration of topological spaces.
\eth

\bpf 
Let $X$ be an object of $\hda$. By Proposition~\ref{ll0} and since the
model category $\top$ is monoidal, one deduces that for any object
${\underline{\alpha}}$ of $\Delta(X^0_{>\alpha})^{op}$, the topological
space $\mathcal{P}_\alpha^-(X)_{\underline{\alpha}}$ is cofibrant.
The pushout product of two cofibrations of topological spaces is
always a cofibration since the model category $\top$ is monoidal.  By
Theorem~\ref{simplification}, it then suffices to prove that for any
loopless object $X$ of $\hda$ such that $(X^0,\leq)$ is locally
finite, for any object $(\alpha_0)$ of $\Delta(X^0_{>\alpha})^{op}$, the
continuous map $L_{(\alpha_0)}\mathcal{P}_\alpha^-(X)\longrightarrow
\mathcal{P}_\alpha^-(X)_{(\alpha_0)}$ is a cofibration of topological spaces. 
Let $X$ be an object of $\hda$. Consider a pushout diagram of flows
with $n\geq 0$ as follows:
\[
\xymatrix{
\glob(\mathbf{S}^{n-1})\fd{}\fr{\phi} & X\fd{}\\
\glob(\mathbf{D}^{n}) \fr{} & \cocartesien Y.}
\] 
One then has to prove that if $X$ satisfies this property, then $Y$
satisfies this property as well. One has $X^0=Y^0$ since the morphism
$\glob(\mathbf{S}^{n-1})\longrightarrow\glob(\mathbf{D}^{n})$
restricts to the identity of $\{\widehat{0},\widehat{1}\}$ on the
$0$-skeletons and since the $0$-skeleton functor $X\mapsto X^0$
preserves colimits~\footnote{One has the canonical bijection $\set(X^0,Z)\iso\dtop(X,T(Z))$ where 
$T(Z)$ is the flow defined by $T(Z)^0=Z$ and for any
$(\alpha,\beta)\in Z\p Z$, $\P_{\alpha,\beta} T(Z)=\{0\}$.}. So one has
the commutative diagram
\[
\xymatrix{
L_{(\alpha_0)}\mathcal{P}_\alpha^-(X)\fr{}\ar@{^(->}[d] & L_{(\alpha_0)}\mathcal{P}_\alpha^-(Y)\fd{}\\
\mathcal{P}_\alpha^-(X)_{(\alpha_0)}\fr{}& \mathcal{P}_\alpha^-(Y)_{(\alpha_0)}}
\]
where the symbol $\xymatrix@1{\ar@{^(->}[r]&}$ means
cofibration. There are two mutually exclusive cases:
\begin{enumerate}
\item $(\phi(\widehat{0}),\phi(\widehat{1}))=(\alpha,\alpha_0)$. 
One then has the situation
\[
\xymatrix{
L_{(\alpha_0)}\mathcal{P}_\alpha^-(X)\ar@{=}[r]\ar@{^(->}[d] & L_{(\alpha_0)}\mathcal{P}_\alpha^-(Y)\fd{}\\
\mathcal{P}_\alpha^-(X)_{(\alpha_0)}\ar@{^{(}->}[r]& \mathcal{P}_\alpha^-(Y)_{(\alpha_0)}}
\]
where the bottom horizontal arrow is a cofibration since it is a
pushout of the morphism of flows
$\glob(\mathbf{S}^{n-1})\longrightarrow\glob(\mathbf{D}^{n})$. So the
continuous map $L_{(\alpha_0)}\mathcal{P}_\alpha^-(Y)\longrightarrow
\mathcal{P}_\alpha^-(Y)_{(\alpha_0)}$ is a cofibration. 
\item $(\phi(\widehat{0}),\phi(\widehat{1}))\neq(\alpha,\alpha_0)$.
  Then, one has the pushout diagram of flows
\[
\xymatrix{
L_{(\alpha_0)}\mathcal{P}_\alpha^-(X)\ar@{->}[r]\ar@{^(->}[d] & L_{(\alpha_0)}\mathcal{P}_\alpha^-(Y)\fd{}\\
\mathcal{P}_\alpha^-(X)_{(\alpha_0)}\ar@{->}[r]& \cocartesien\mathcal{P}_\alpha^-(Y)_{(\alpha_0)}}
\]
So the continuous map $L_{(\alpha_0)}\mathcal{P}_\alpha^-(Y)\longrightarrow
\mathcal{P}_\alpha^-(Y)_{(\alpha_0)}$ is again a cofibration. In this situation, 
it may happen that
$L_{(\alpha_0)}\mathcal{P}_\alpha^-(X)=L_{(\alpha_0)}\mathcal{P}_\alpha^-(Y)$.
\end{enumerate}
The proof is complete with Proposition~\ref{compgen}, and because the
canonical morphism of flows $X^0 \longrightarrow X$ is a
relative $I^{gl}$-cell complex, and at last because the property above
is clearly satisfied for $X=X^0$.
\epf

\section{The end of the proof}
\label{end}

The two following classical results about classifying spaces are going
to be very useful.

\bp (for example \cite{ref_model2} Proposition~18.1.6) \label{calculconcret} 
Let $\C$ be a small category. Then the homotopy colimit of the
terminal object of $\top^\C$ is homotopy equivalent to the classifying
space of the opposite of the indexing category $\C$. 
\ep

\bp (for example \cite{ref_model2} Proposition~14.3.13) \label{calculconcret2}
Let $\C$ be a small category having a terminal object. Then the classifying 
space of $\C$ is contractible. 
\ep

\bth \label{finpm} Let $X$ be a loopless object of $\hda$ such that
$(X^0,\leq)$ is locally finite. Assume there is an element
$\widehat{1}$ such that $\alpha\leq \widehat{1}$ for all $\alpha\in
X^0$. For any $1$-simplex $(\alpha,\beta)$ of $\Delta(X^0)$, let us
suppose that $\P_{\alpha,\beta}X$ is weakly contractible. Then
$\hop^-_\alpha X$ has the homotopy type of a point for any $\alpha\in
X^0\backslash \{\widehat{1}\}$.  \eth

Of course, with the hypothesis of the theorem, the topological space
$\hop^-_{\widehat{1}} X$ is the empty space.

\bpf One has the sequence of weak homotopy equivalences (where $Q$ is the cofibrant replacement functor, where $\mathbf{1}$ is the terminal diagram and with $\alpha\in X^0\backslash
\{\widehat{1}\}$)
{\small 
\begin{align*}
&\hop^-_\alpha X& \\
&\simeq\P_\alpha^-Q(X) & \hbox{ by definition of the homotopy branching space}\\
&\simeq \liminj_{\Delta(X^0_{>\alpha})^{op}} \mathcal{P}_\alpha^-(Q(X)) & \hbox{ by Theorem~\ref{calculPm}}\\
&\simeq \holiminj_{\Delta(X^0_{>\alpha})^{op}} \mathcal{P}_\alpha^-(Q(X)) & \hbox{ by Theorem~\ref{equi} and Theorem~\ref{cascofibrant}}\\
&\simeq \holiminj_{\Delta(X^0_{>\alpha})^{op}} \mathbf{1} & \hbox{ by homotopy invariance of the homotopy colimit}\\
&\simeq B\lp\Delta(X^0_{>\alpha})\rp & \hbox{ by Proposition~\ref{calculconcret}}\\
&\simeq B(X^0_{>\alpha}) & \hbox{ since the barycentric subdivision is a homotopy invariant}\\
&\simeq B(]\alpha,\widehat{1}]) & \hbox{ since $X^0$ has exactly one maximal point}\\
&\simeq \{0\} & \hbox{ by Proposition~\ref{calculconcret2}.}
\end{align*}}
\epf

\bp \label{universalmp}(\cite{exbranch} Proposition~A.1)
Let $X$ be a flow. There exists a topological space $\P^+X$ unique up
to homeomorphism and a continuous map $h^+:\P X\longrightarrow \P^+ X$
satisfying the following universal property:
\begin{enumerate}
\item For any $x$ and $y$ in $\P X$ such that $t(x)=s(y)$, the equality
$h^+(y)=h^+(x*y)$ holds.
\item Let $\phi:\P X\longrightarrow Y$ be a
continuous map such that for any $x$ and $y$ of $\P X$ such that
$t(x)=s(y)$, the equality $\phi(y)=\phi(x*y)$ holds. Then there exists a
unique continuous map $\overline{\phi}:\P^+X\longrightarrow Y$ such that
$\phi=\overline{\phi}\circ h^+$.
\end{enumerate}
Moreover, one has the homeomorphism
\[\P^+X\iso \bigsqcup_{\alpha\in X^0} \P^+_\alpha X\]
where $\P^+_\alpha X:=h^+\lp \bigsqcup_{\beta\in
X^0} \P^+_{\alpha,\beta} X\rp$. The mapping $X\mapsto \P^+X$
yields a functor $\P^+$ from $\dtop$ to $\top$. 
\ep

Roughly speaking, the merging space of a flow is the space of germs
of non-constant execution paths ending in the same way.

\bd \label{defplus}
Let $X$ be a flow. The topological space $\P^+X$ is called the {\rm
merging space} of the flow $X$. The functor $\P^+$ is called the 
{\rm merging space functor}. \ed

\bth  (\cite{exbranch} Theorem~A.4)
The merging space functor $\P^+:\dtop\longrightarrow \top$ is a left
Quillen functor.
\eth

\bd
The {\rm homotopy merging space} $\hop^+ X$ of a flow $X$ is by
definition the topological space $\P^+Q(X)$.  If $\alpha\in X^0$, let
$\hop^+_\alpha X=\P^+_\alpha Q(X)$.
\ed

\bth 
\label{resultat1} 
Let $\vec{D}$ be a full directed ball with initial state $\widehat{0}$ and
final state $\widehat{1}$.  Then one has the homotopy equivalences
\begin{enumerate}
\item $\hop^-_\alpha \vec{D}\simeq \{0\}$ for any $\alpha\in \vec{D}^0\backslash \{\widehat{1}\}$ and 
$\hop^-_{\widehat{1}} \vec{D}=\varnothing$
\item $\hop^+_\alpha \vec{D}\simeq \{0\}$ for any $\alpha\in \vec{D}^0\backslash \{\widehat{0}\}$ and 
$\hop^+_{\widehat{0}} \vec{D}=\varnothing$. 
\end{enumerate}
\eth

\bpf The equalities $\hop^-_{\widehat{1}} \vec{D}=\varnothing$ and $\hop^+_{\widehat{0}} \vec{D}=\varnothing$ are 
obvious. The homotopy equivalence $\hop^-_\alpha \vec{D}\simeq \{0\}$ for
any $\alpha\in \vec{D}^0\backslash \{{\widehat{1}}\}$ is the result of
Theorem~\ref{finpm}. Let us consider the opposite flow $\vec{D}^{op}$ of $\vec{D}$
defined as follows:
\begin{enumerate}
\item $(\vec{D}^{op})^0=\vec{D}^0$
\item $\P_{\alpha,\beta}\vec{D}^{op}=\P_{\beta,\alpha}\vec{D}$ with $t^{op}(\gamma)=s(\gamma)$ and 
$s^{op}(\gamma)=t(\gamma)$. 
\end{enumerate}
The weak S-homotopy equivalence $Q(\vec{D})\longrightarrow \vec{D}$ from the
cofibrant replacement of $\vec{D}$ to $\vec{D}$ becomes a weak S-homotopy
equivalence $Q(\vec{D})^{op}\longrightarrow \vec{D}^{op}$. Since one has the
isomorphism $\glob(Z)^{op}\iso \glob(Z)$ for any topological space $Z$
(in particular, for $Z=\mathbf{S}^{n-1}$ and $Z=\mathbf{D}^n$ for
all $n\geq 0$), then the transfinite composition
$\varnothing\longrightarrow Q(\vec{D})$ of pushouts of morphisms of
\[\{\glob(\mathbf{S}^{n-1})\longrightarrow \glob(\mathbf{D}^n),n\geq 0\}\cup \{R,C\}\] 
allows to view $\varnothing\longrightarrow Q(\vec{D})^{op}$ as the
transfinite composition of pushouts of the same set of
morphisms. Therefore, the flow $Q(\vec{D})^{op}$ is a cofibrant replacement
functor of $\vec{D}^{op}$. So one has the homotopy equivalences 
\[\hop^+_\alpha \vec{D}\simeq \P^+_\alpha Q(\vec{D})\simeq \P^-_\alpha Q(\vec{D})^{op} 
\simeq \P^-_\alpha Q(\vec{D}^{op})\simeq \hop^-_\alpha \vec{D}^{op}.\] 
Thus, if $\alpha$ is not the final state of $\vec{D}^{op}$, that is the
initial state of $\vec{D}$, then we are reduced to verifying that $\vec{D}^{op}$
is a full directed ball as well. The latter fact is clear.
\epf

\section{The branching and merging homologies of a flow}
\label{rappelHMP}

We recall in this section the definition of the branching and merging
homologies.

\bd\label{hombrdef} \cite{exbranch}
Let $X$ be a flow. Then the $(n+1)$-st branching homology group
$H_{n+1}^-(X)$ is defined as the $n$-th homology group of the
augmented simplicial set $\mathcal{N}^-_*(X)$ defined as follows:
\begin{enumerate}
\item $\mathcal{N}^-_n(X)=\sing_n(\hop^-X)$ for $n\geq 0$
\item $\mathcal{N}^-_{-1}(X)=X^0$
\item the augmentation map $\epsilon:\sing_0(\hop^-X)\longrightarrow X^0$
is induced by the mapping $\gamma\mapsto s(\gamma)$ from $\hop^-X=\sing_0(\hop^-X)$
to $X^0$
\end{enumerate}
where $\sing(Z)$ denotes the singular simplicial nerve of a given
topological space $Z$ \cite{MR2001d:55012}. In other terms, 
\begin{enumerate}
\item for $n\geq 1$, $H_{n+1}^-(X):=H_n(\hop^-X)$
\item  $H_1^-(X):=\ker(\epsilon)/\im\lp\partial:\mathcal{N}^-_1(X)\rightarrow
\mathcal{N}^-_0(X)\rp$
\item $H_0^-(X):=\Z(X^0)/\im(\epsilon)$.
\end{enumerate}
where $\partial$ is the simplicial differential map, where $\ker(f)$ is the kernel 
of $f$ and where $\im(f)$ is the image of $f$. 
\ed

For any flow $X$, $H_0^-(X)$ is the free abelian group generated by
the final states of $X$.

\bd \cite{exbranch}
Let $X$ be a flow. Then the $(n+1)$-st merging homology group
$H_{n+1}^+(X)$ is defined as the $n$-th homology group of the
augmented simplicial set $\mathcal{N}^+_*(X)$ defined as follows:
\begin{enumerate}
\item $\mathcal{N}^+_n(X)=\sing_n(\hop^+X)$ for $n\geq 0$
\item $\mathcal{N}^+_{-1}(X)=X^0$
\item the augmentation map $\epsilon:\sing_0(\hop^+X)\longrightarrow X^0$
is induced by the mapping $\gamma\mapsto t(\gamma)$ from $\hop^+X=\sing_0(\hop^+X)$
to $X^0$
\end{enumerate}
where $\sing(Z)$ denotes the singular simplicial nerve of a given topological space 
$Z$. In other terms, 
\begin{enumerate}
\item for $n\geq 1$, $H_{n+1}^+(X):=H_n(\hop^+X)$
\item  $H_1^+(X):=\ker(\epsilon)/\im\lp\partial:\mathcal{N}^+_1(X)\rightarrow
\mathcal{N}^+_0(X)\rp$
\item $H_0^+(X):=\Z(X^0)/\im(\epsilon)$.
\end{enumerate}
where $\partial$ is the simplicial differential map, where $\ker(f)$ is the kernel 
of $f$ and where $\im(f)$ is the image of $f$. 
\ed

For any flow $X$, $H_0^+(X)$ is the free abelian group generated by the initial
states of $X$.

\section{Preservation of the branching and merging homologies}
\label{preHMP}

\bd \cite{model2}
Let $X$ be a flow. Let $A$ and $B$ be two subsets of $X^0$. One says that
$A$ is {\rm surrounded} by $B$ (in $X$) if for any $\alpha\in A$,
either $\alpha \in B$ or there exists execution paths $\gamma_1$ and
$\gamma_2$ of $\P X$ such that $s(\gamma_1)\in B$,
$t(\gamma_1)=s(\gamma_2)=\alpha$ and $t(\gamma_2)\in B$. We denote
this situation by $A\lll B$.  \ed

\bth \label{pre1} 
Let $f:X\longrightarrow Y$ be a generalized T-homotopy
equivalence. Then the morphism of flows $f$ satisfies the 
following conditions (with $\epsilon=\pm$): 
\begin{enumerate}
\item $Y^0\lll f(X^0)$
\item for any $\alpha\in X^0$, $f$ induces a weak homotopy equivalence  
$\hop^\epsilon_\alpha X\simeq \hop^\epsilon_{f(\alpha)} Y$
\item for any $\alpha\in Y^0\backslash f(X^0)$, the topological space 
$\hop^\epsilon_{\alpha} Y$ is contractible. 
\end{enumerate}
\eth

\bpf First of all, let us suppose that $f$ is a pushout of the form
\[
\xymatrix{
Q(F(P_1))\ar@{^{(}->}[d]_{Q(F(u))}\fr{}& X\fd{f}\\ Q(F(P_2)) \fr{} & Y
\cocartesien}
\]
where $P_1$ and $P_2$ are two finite bounded posets and where
$u:P_1\longrightarrow P_2$ belongs to $\mathcal{T}$
(Definition~\ref{definitiondeT}). Let us factor the morphism of flows
$Q(F(P_1))\longrightarrow X$ as a composite of a cofibration
$Q(F(P_1))\longrightarrow W$ followed by a trivial fibration
$W\longrightarrow X$.  Then one obtains the commutative diagram of
flows
\[
\xymatrix{
Q(F(P_1))\ar@{^{(}->}[d]_{Q(F(u))}\ar@{^{(}->}[r] & W \ar@{^{(}->}[d] \ar@{->>}[r]^{\simeq}& X \fd{f}\\
Q(F(P_2)) \ar@{^{(}->}[r] & \cocartesien T \fr{\simeq} & \cocartesien Y.}
\]
The morphism $T\longrightarrow Y$ of the diagram above is a weak
S-homotopy equivalence since the model category $\dtop$ is left proper
by \cite{2eme} Theorem~6.4. So the flows $W$ and $X$ (resp. $T$ and
$Y$) have same homotopy branching and merging spaces and we are
reduced to the following situation:
\[
\xymatrix{
Q(F(P_1))\ar@{^{(}->}[d]_{Q(F(u))}\ar@{^{(}->}[r]& X\ar@{^{(}->}[d]_{f}\\
Q(F(P_2)) \ar@{^{(}->}[r] & Y \cocartesien}
\]
where the square is both a pushout and a homotopy pushout diagram of
flows. The $0$-skeleton functor gives rise to the commutative diagram 
of set maps: 
\[
\xymatrix{
P_1\ar@{->}[d]_{u}\fr{v}& X^0\ar@{->}[d]\\
P_2 \fr{w} & Y^0. \cocartesien}
\]
Thus, one obtains the commutative diagram of topological spaces
($\epsilon\in\{-1,+1\}$)
\[
\xymatrix{
\bigsqcup_{\beta\in v^{-1}(\alpha)}\P^\epsilon_\beta Q(F(P_1))\ar@{^{(}->}[d]^{\simeq}\ar@{^{(}->}[r]& \P^\epsilon_\alpha X\ar@{^{(}->}[d]^{\simeq}\\
\bigsqcup_{\beta\in w^{-1}(f(\alpha))}\P^\epsilon_\beta Q(F(P_2)) \ar@{^{(}->}[r] & \P^\epsilon_{f(\alpha)} Y. \cocartesien}
\]
The left vertical arrow is a weak homotopy equivalence for the following reasons: 
\begin{enumerate}
\item 
Theorem~\ref{resultat1} says that each component of the domain 
and of the codomain is weakly contractible, or empty. And since 
$u(\widehat{0})=\widehat{0}$ and $u(\widehat{1})=\widehat{1}$, a component 
$\P^\epsilon_\beta Q(F(P_1))$ is empty (resp. weakly contractible) if and only 
if $\P^\epsilon_{u(\beta)} Q(F(P_2))$ is empty (resp. weakly contractible). 
\item The map $u$ is one-to-one and therefore, the restriction 
\[u: v^{-1}(\alpha) \longrightarrow w^{-1}(f(\alpha))\] is bijective. 
\end{enumerate}
The left vertical arrow is also a cofibration. So the right vertical
arrow $\P^\epsilon_\alpha X\longrightarrow \P^\epsilon_{f(\alpha)} Y$
is a trivial cofibration as well since the functors $\P^\epsilon$ with
$\epsilon\in\{-1,+1\}$ are both left Quillen functors.

Let $\alpha\in Y^0\backslash f(X^0)$, that is to say $\alpha\in
P_2\backslash u(P_1)$. Then one obtains the pushout diagram of
topological spaces
\[
\xymatrix{
\varnothing=\P^\epsilon_\alpha Q(\vI)\ar@{^{(}->}[d]_{}\ar@{^{(}->}[r]& \varnothing=\P^\epsilon_\alpha X\ar@{^{(}->}[d]\\
\P^\epsilon_\alpha Q(D) \ar@{^{(}->}[r] & \P^\epsilon_\alpha Y. \cocartesien}
\]
So by Theorem~\ref{resultat1} again, one deduces that
$\P^\epsilon_\alpha Y$ is contractible as soon as $\alpha\in Y^0\backslash f(X^0)$.

Now let us suppose that $f:X\longrightarrow Y$ is a transfinite
composition of morphisms as above. Then there exists an ordinal
$\lambda$ and a $\lambda$-sequence $Z:\lambda\longrightarrow \dtop$
with $Z_0=X$, $Z_\lambda=Y$ and the morphism $Z_0\longrightarrow
Z_\lambda$ is equal to $f$. Since for any $u\in\mathcal{T}$, the
morphism of flows $Q(F(u))$ is a cofibration, the morphism
$Z_\mu\longrightarrow Z_{\mu+1}$ is a cofibration for any
$\mu<\lambda$.  Since the model category $\dtop$ is left proper by
\cite{2eme} Theorem~6.4, there exists by \cite{ref_model2}
Proposition~17.9.4 a $\lambda$-sequence
$\widetilde{Z}:\lambda\longrightarrow \dtop$ and a morphism of
$\lambda$-sequences $\widetilde{Z}\longrightarrow Z$ such that for any
$\mu\leq \lambda$, the flow $\widetilde{Z}_\mu$ is cofibrant and the
morphism $\widetilde{Z}_\mu\longrightarrow Z_\mu$ is a weak S-homotopy
equivalence. So for any $\mu\leq \lambda$, one has $\P^\epsilon
\widetilde{Z}_\mu\simeq \hop^\epsilon Z_\mu$ and for any
$\mu<\lambda$, the continuous map $\P^\epsilon
\widetilde{Z}_\mu\longrightarrow \P^\epsilon \widetilde{Z}_{\mu+1}$ is
a cofibration. So for a given $\alpha\in Z_0^0=X^0$, the continuous
map $\hop^\epsilon_\alpha X\longrightarrow \hop^\epsilon_{f(\alpha)}Y$
is a transfinite composition of trivial cofibrations, and therefore a
trivial cofibration as well.

The same argument proves that the continuous map $\hop^\epsilon_\alpha
Z_\mu\longrightarrow \hop^\epsilon_{\alpha'}Y$ is a trivial
cofibration for any $\mu\leq \lambda$ where $\alpha'\in Y^0$ is the
image of $\alpha\in Z_\mu^0$ by the morphism $Z_\mu\longrightarrow
Y$. Let $\alpha\in Y^0\backslash f(X^0)$. Consider the set of ordinals
\[\left\{\mu\leq\lambda, \exists \beta_\mu\in Z_\mu \hbox{ mapped to }\alpha\right\}.\]
This non-empty set (it contains at least $\lambda$) has a smallest
element $\mu_0$. The ordinal $\mu_0$ cannot be a limit ordinal.
Otherwise, one would have $Z_{\mu_0}=\liminj_{\mu<\mu_0} Z_{\mu}$ and
therefore there would exist a $\beta_\mu$ mapped to $\beta_{\mu_0}$
for some $\mu<\mu_0$: contradiction. So one can write $\mu_0=\mu_1+1$.
Then $\hop^\epsilon_{\beta_{\mu_0}}Z_{\mu_0}$ is contractible because
of the first part of the proof applied to the morphism
$Z_{\mu_1}\longrightarrow Z_{\mu_0}$. Therefore,
$\hop^\epsilon_{\alpha'}Y$ is contractible as well.

The condition $Y^0\lll f(X^0)$ is always clearly satisfied.

It remains the case where $f$ is a retract of a generalized
T-equivalence of the preceding kinds. The result follows from the fact
that everything is functorial and that the retract of a weak homotopy
equivalence (resp. a non-empty set) is a weak homotopy equivalence
(resp. a non-empty set).
\epf

\begin{cor} \label{enfinfin}
Let $f:X\longrightarrow Y$ be a generalized T-homotopy
equivalence. Then for any $n\geq 0$, one has the isomorphisms
$H_n^-(f):H_n^-(X)\longrightarrow H_n^-(Y)$ and
$H_n^+(f):H_n^+(X)\longrightarrow H_n^+(Y)$.
\end{cor}

\bpf 
This is the same proof as for \cite{exbranch} Proposition~7.4 (the
word contractible being replaced by singleton).
\epf

\section{Conclusion}

This new definition of T-homotopy equivalence seems to be well-behaved
because it preserves the branching and merging homology theories. For
an application of this new approach of T-homotopy, see the proof of an
analogue of Whitehead's theorem for the full dihomotopy relation in
\cite{hocont}.

\appendix

\section{Elementary remarks about flows}
\label{limelm}

\bp \label{precision}
(\cite{model3} Proposition~15.1) If one has the pushout of flows
\[
\xymatrix{
\glob(\de Z) \fr{\phi}\fd{} & A \fd{} \\
\glob(Z) \fr{} & \cocartesien M }
\]
then the continuous map $\P A\longrightarrow \P M$ is a transfinite
composition of pushouts of continuous maps of the form
$\id\p\dots\p\id\p f \p\id \p\dots\p \id$ where
$f:\P_{\phi(\widehat{0}),\phi(\widehat{1})} A\longrightarrow T$ is the
canonical inclusion obtained with the pushout diagram of topological
spaces
\[
\xymatrix{
\de Z\fr{}\fd{} & \P_{\phi(\widehat{0}),\phi(\widehat{1})}A\fd{} \\
Z \fr{} & \cocartesien T. }
\]
\ep

\bp\label{ll05}
Let $Y$ be a flow such that $\P Y$ is a cofibrant topological
space. Let $f:Y\longrightarrow Z$ be a pushout of a morphism of
$I^{gl}_+$. Then the topological space $\P Z$ is cofibrant.
\ep

\bpf By hypothesis, $f$ is the pushout of a morphism of flows
$g\in I^{gl}_+$. So one has the pushout of flows
\[
\xymatrix{
A\fd{g} \fr{\phi} & Y\ar@{->}[d]^-{f}\\
B \ar@{->}[r]_-{\psi} & Z. \cocartesien
}
\]
If $f$ is a pushout of $C:\varnothing\subset \{0\}$, then $\P Z=\P Y$.
Therefore, the space $\P Z$ is cofibrant. If $f$ is a pushout of
$R:\{0,1\}\rightarrow \{0\}$ and if $\phi(0)=\phi(1)$, then $\P Z=\P
Y$ again.  Therefore, the space $\P Z$ is also cofibrant. If $f$ is a
pushout of $R:\{0,1\}\rightarrow \{0\}$ and if $\phi$ is one-to-one,
then one has the homeomorphism \beas \P Z&\iso&\P Y\sqcup
\bigsqcup_{r\geq 0} \lp
\P_{.,\phi(0)}Y\p\P_{\phi(1),\phi(0)}Y\p \P_{\phi(1),\phi(0)}Y\p \dots\hbox{($r$ times)}  \p \P_{\phi(1),.}Y  \rp \\
&&\sqcup \bigsqcup_{r\geq 0} \lp
\P_{.,\phi(1)}Y\p\P_{\phi(0),\phi(1)}Y\p \P_{\phi(0),\phi(1)}Y\p
\dots\hbox{($r$ times)} \p \P_{\phi(0),.}Y \rp.  \eeas Therefore, the
space $\P Z$ is again cofibrant since the model category $\top$ is
monoidal.  It remains the case where $g$ is the inclusion
$\glob(\mathbf{S}^{n-1})\subset \glob(\mathbf{D}^n)$ for some $n\geq
0$. Consider the pushout of topological spaces
\[
\xymatrix{
\mathbf{S}^{n-1}\fd{g} \fr{\P \phi} & \P_{\phi(\widehat{0}),\phi(\widehat{1})}Y\ar@{->}[d]^-{f}\\
\mathbf{D}^n \ar@{->}[r]_-{\P \psi} & T. \cocartesien
}
\]
By Proposition~\ref{precision}, the continuous map $\P
Y\longrightarrow \P Z$ is a transfinite composition of pushouts of
continuous maps of the form $\id\p\id\p\dots \p f \p \dots \p
\id\p\id$ where $f$ is a cofibration and the identities maps are the
identity maps of cofibrant topological spaces. So it suffices to
notice that if $k$ is a cofibration and if $X$ is a cofibrant
topological space, then $\id_X\p k$ is still a cofibration since the
model category $\top$ is monoidal.  \epf

\bp\label{ll0} 
Let $X$ be a cofibrant flow.  Then for any $(\alpha,\beta)\in X^0\p
X^0$, the topological space $\P_{\alpha,\beta}X$ is cofibrant. 
\ep

\bpf 
A cofibrant flow $X$ is a retract of a $I^{gl}_+$-cell complex $Y$ and
$\P X$ becomes a retract of $\P Y$. So it suffices to show that $\P Y$
is cofibrant. Proposition~\ref{ll05} completes the proof.
\epf

\section{Calculating pushout products}
\label{calpush}

\begin{lem} \label{colimproduit} 
Let $D:I\longrightarrow \top$ and $E:J\longrightarrow \top$ be two
diagrams in a complete cocomplete cartesian closed category. Let $D\p
E:I\p J:\longrightarrow \top$ be the diagram of topological spaces
defined by $(D\p E)(x,y):=D(x)\p E(y)$ if $(x,y)$ is either an object
or an arrow of the small category $I\p J$. Then one has $\liminj (D\p
E)\iso(\liminj D)\p (\liminj E)$.
\end{lem}

\bpf One has $\liminj (D\p E)\iso \liminj_i(\liminj_j D(i)\p E(j))$ by 
\cite{MR1712872}. And  one has $\liminj_j (D(i)\p E(j))\iso D(i)\p (\liminj E)$ 
since the category is cartesian closed. So $\liminj (D\p E)\iso
\liminj_i (D(i)\p (\liminj E))\iso (\liminj D)\p (\liminj E)$. \epf

\begin{nota} \label{notapush}
If $f:U\longrightarrow V$ and $g:W\longrightarrow X$ are two morphisms
of a complete cocomplete category, then let us denote by $f\square g:
(U\p X) \sqcup_{(U\p W)} (V\p W)\longrightarrow V\p X$ the {\rm
pushout product} of $f$ and $g$. The notation $f_0\square\dots \square
f_p$ is defined by induction on $p$ by $f_0\square\dots \square
f_p:=(f_0\square\dots \square f_{p-1})\square f_p$.
\end{nota}

\bth\label{calculpushout}  (Calculating a pushout product of several morphisms)
Let $f_i:A_i\longrightarrow B_i$ for $0\leq i\leq p$ be $p+1$ morphisms
of a complete cocomplete cartesian closed category $\C$. Let $S\subset
\{0,\dots,p\}$. Let 
\[C_p(S):=\lp \prod_{i\in S} B_i\rp \p \lp \prod_{i\notin S} A_i\rp.\]
If $S$ and $T$ are two subsets of $\{0,\dots,p\}$ such that $S\subset
T$, let $C_p(i_S^T):C_p(S)\longrightarrow C_p(T)$ be the morphism 
\[\lp\prod_{i\in S}\id_{B_i}\rp\p \lp\prod_{i\in T\backslash S} f_i\rp\p 
\lp\prod_{i\notin T} \id_{A_i}\rp. \]
Then:  
\begin{enumerate} 
\item the mappings $S\mapsto C_p(S)$ and $i_S^T\mapsto C_p(i_S^T)$ 
give rise to a functor from $\Delta(\{0,\dots,p\})$ (the order
complex of the poset $\{0,\dots,p\}$) to $\C$
\item there exists a canonical morphism 
{
\[\liminj_{S\subsetneqq \{0,\dots,p\}} 
C_p(S)\longrightarrow C_p(\{0,\dots,p\}).\]}
and it is equal to the morphism $f_0\square\dots \square f_p$. 
\end{enumerate}
\eth

\bpf The first assertion is clear. Moreover, for any subset $S$ and $T$ of 
$\{0,\dots,p\}$ such that $S\subset T$, the diagram
\[
\xymatrix{
S\fr{}\fd{} & \{0,\dots,p\} \\
T\ar@{->}[ru]&}
\]
is commutative since there is at most one morphism between two objects
of the category $\Delta(\{0,\dots,p\})$, hence the existence of the
morphism {\[\liminj_{S\subsetneqq \{0,\dots,p\}} C_p(S)\longrightarrow
  C(\{0,\dots,p\}).\]}

The second assertion is clear for $p=0$ and $p=1$. We are going to
prove it by induction on $p$. By definition, the morphism
$f_0\square\dots \square f_{p+1}$ is the canonical morphism from
{\small\[\lp\lp\liminj_{S\subsetneqq \{0,\dots,p\}} C_p(S)\rp \p B_{p+1}\rp
  \sqcup_{\lp\lp\liminj_{S\subsetneqq \{0,\dots,p\}} C_p(S)\rp \p
    A_{p+1}\rp} \lp C_p(\{0,\dots,p\})\p A_{p+1}\rp\]} to
$B_0\p\dots\p B_{p+1}$.  By Lemma~\ref{colimproduit}, the source of
the morphism $f_0\square\dots \square f_{p+1}$ is then equal to
\[\lp\liminj_{p+1\in S\subsetneqq \{0,\dots,p+1\}} 
C_{p+1}(S)\rp \sqcup_{\lp\liminj_{S\subsetneqq \{0,\dots,p\}}
  C_{p+1}(S)\rp} \lp C_{p+1}(\{0,\dots,p\})\rp,\] an the latter is
equal to $\liminj_{S\subsetneqq \{0,\dots,p+1\}} C_{p+1}(S)$.
\epf

\section{Mixed transfinite composition of pushouts and cofibrations}
\label{mixcomp}

\bp \label{compgen} Let $\mathcal{M}$ be a model category. Let
$\lambda$ be an ordinal.  Let $(f_\mu:A_\mu\longrightarrow
B_\mu)_{\mu<\lambda}$ be a $\lambda$-sequence of morphisms of
$\mathcal{M}$. For $\mu<\lambda$, suppose that $A_\mu\rightarrow
A_{\mu+1}$ is an isomorphism or the diagram of objects of
$\mathcal{M}$
\[
\xymatrix{
A_\mu \fr{} \fd{f_\mu} & A_{\mu+1}\fd{}\\
B_\mu \ar@{->}[r]  & B_{\mu+1}
}
\]
is a pushout, and for $\mu<\lambda$ suppose also that the map
$B_\mu\longrightarrow B_{\mu+1}$ is a cofibration. Then: if
$f_0:A_0\longrightarrow B_0$ is a cofibration, then
$f_\lambda:A_\lambda\longrightarrow B_\lambda$ is a cofibration as
well, where of course $A_\lambda:=\liminj A_\mu$ and
$B_\lambda:=\liminj B_\mu$.  \ep

\bpf It is clear that if $f_\mu:A_\mu\longrightarrow B_\mu$ is a cofibration, then 
$f_{\mu+1}:A_{\mu+1}\longrightarrow B_{\mu+1}$ is a cofibration as
well. It then suffices to prove that if $\nu\leq\lambda$ is a limit
ordinal such that $f_\mu:A_\mu\longrightarrow B_\mu$ is a cofibration
for any $\mu<\nu$, then $f_\nu:A_\nu\longrightarrow B_\nu$ is a
cofibration as well.  Consider a commutative diagram
\[
\xymatrix{
A_\nu \fr{} \fd{f_\nu} & C\fd{}\\
B_\nu \fr{}  \ar@{-->}[ru]^{k}& D
}
\]
where $C\longrightarrow D$ is a trivial fibration of
$\mathcal{M}$. Then one has to find $k:B_\nu\longrightarrow C$ making
both triangles commutative. Recall that by hypothesis,
$f_\nu=\liminj_{\mu<\nu} f_\mu$. Since $f_0$ is a cofibration, there 
exists a map $k_0$ making both triangles of the diagram 
\[
\xymatrix{
A_0 \fr{} \fd{f_0} & C\fd{}\\
B_0 \fr{}  \ar@{-->}[ru]^{k_0}& D
}
\]
commutative. Let us suppose $k_\mu$ constructed. There are two
cases. Either the diagram
\[
\xymatrix{
A_\mu \fr{} \fd{f_\mu} & A_{\mu+1}\fd{}\\
B_\mu \ar@{^(->}[r]  & B_{\mu+1}
}
\] 
is a pushout, and one can construct a morphism $k_{\mu+1}$ making both
triangles of the diagram
\[
\xymatrix{
A_{\mu+1} \fr{} \fd{f_{\mu+1}} & C\fd{}\\
B_{\mu+1} \fr{}  \ar@{-->}[ru]^{k_{\mu+1}}& D
}
\]
commutative and such that the composite $B_\mu \longrightarrow
B_{\mu+1} \longrightarrow C$ is equal to $k_\mu$ by using the
universal property satisfied by the pushout. Or the morphism
$A_\mu\rightarrow A_{\mu+1}$ is an isomorphism. In that latter case,
consider the commutative diagram
\[
\xymatrix{
B_\mu \fr{} \ar@{^(->}[d] \fr{k_\mu} & C\fd{}\\
B_{\mu+1} \fr{}  & D
}
\] 
Since the morphism $B_\mu \longrightarrow B_{\mu+1}$ is a cofibration,
there exists $k_{\mu+1} : B_{\mu+1} \longrightarrow C$ making the two
triangles of the latter diagram commutative. So, once again, the
composite $B_\mu \longrightarrow B_{\mu+1} \longrightarrow C$ is equal
to $k_\mu$.

The map $k:=\liminj_{\mu<\nu} k_\mu$ is a solution.
\epf


\begin{thebibliography}{Gau06b}

\bibitem[Bro88]{MR90k:54001}
R.~Brown.
\newblock {\em Topology}.
\newblock Ellis Horwood Ltd., Chichester, second edition, 1988.
\newblock A geometric account of general topology, homotopy types and the
  fundamental groupoid.

\bibitem[DS95]{MR1361887}
W.~G. Dwyer and J.~Spali{\'n}ski.
\newblock Homotopy theories and model categories.
\newblock In {\em Handbook of algebraic topology}, pages 73--126.
  North-Holland, Amsterdam, 1995.

\bibitem[Gau00]{ConcuToAlgTopo}
P.~Gaucher.
\newblock From concurrency to algebraic topology.
\newblock In {\em Electronic Notes in Theoretical Computer Science}, volume~39,
  page 19pp, 2000.

\bibitem[Gau03]{model3}
P.~Gaucher.
\newblock A model category for the homotopy theory of concurrency.
\newblock {\em Homology Homotopy Appl.}, 5(1):549--599 (electronic), 2003.

\bibitem[Gau05a]{model2}
P.~Gaucher.
\newblock Comparing globular complex and flow.
\newblock {\em New York J. Math.}, 11:97--150 (electronic), 2005.

\bibitem[Gau05b]{exbranch}
P.~Gaucher.
\newblock Homological properties of non-deterministic branchings of mergings in
  higher dimensional automata.
\newblock {\em Homology Homotopy Appl.}, 7(1):51--76 (electronic), 2005.

\bibitem[Gau05c]{1eme}
P.~Gaucher.
\newblock T-homotopy and refinement of observation ({I}) : {I}ntroduction.
\newblock preprint ArXiv math.AT, to appear in ENTCS, 2005.

\bibitem[Gau05d]{2eme}
P.~Gaucher.
\newblock T-homotopy and refinement of observation ({II}) : {A}dding new
  {T}-homotopy equivalences.
\newblock preprint ArXiv math.AT, 2005.

\bibitem[Gau06a]{hocont}
P.~Gaucher.
\newblock Inverting weak dihomotopy equivalence using homotopy continuous flow.
\newblock {\em Theory Appl. Categ.}, 16:No. 3, 59--83 (electronic), 2006.

\bibitem[Gau06b]{4eme}
P.~Gaucher.
\newblock T-homotopy and refinement of observation ({IV}) : Invariance of the
  underlying homotopy type.
\newblock {\em New York J. Math.}, 12:63--95 (electronic), 2006.

\bibitem[GG03]{diCW}
P.~Gaucher and E.~Goubault.
\newblock Topological deformation of higher dimensional automata.
\newblock {\em Homology Homotopy Appl.}, 5(2):39--82 (electronic), 2003.
\newblock Algebraic topological methods in computer science (Stanford, CA,
  2001).

\bibitem[GJ99]{MR2001d:55012}
P.~G. Goerss and J.~F. Jardine.
\newblock {\em Simplicial homotopy theory}.
\newblock Birkh\"auser Verlag, Basel, 1999.

\bibitem[Gou03]{survol}
E.~Goubault.
\newblock Some geometric perspectives in concurrency theory.
\newblock {\em Homology Homotopy Appl.}, 5(2):95--136 (electronic), 2003.
\newblock Algebraic topological methods in computer science (Stanford, CA,
  2001).

\bibitem[Hir03]{ref_model2}
P.~S. Hirschhorn.
\newblock {\em Model categories and their localizations}, volume~99 of {\em
  Mathematical Surveys and Monographs}.
\newblock American Mathematical Society, Providence, RI, 2003.

\bibitem[Hov99]{MR99h:55031}
M.~Hovey.
\newblock {\em Model categories}.
\newblock American Mathematical Society, Providence, RI, 1999.

\bibitem[Lew78]{Ref_wH}
L.~G. Lewis.
\newblock {\em The stable category and generalized Thom spectra}.
\newblock PhD thesis, University of Chicago, 1978.

\bibitem[May99]{MR2000h:55002}
J.~P. May.
\newblock {\em A concise course in algebraic topology}.
\newblock University of Chicago Press, Chicago, IL, 1999.

\bibitem[ML98]{MR1712872}
S.~Mac~Lane.
\newblock {\em Categories for the working mathematician}.
\newblock Springer-Verlag, New York, second edition, 1998.

\bibitem[Qui67]{MR36:6480}
D.~Quillen.
\newblock {\em Homotopical algebra}.
\newblock Springer-Verlag, Berlin, 1967.

\bibitem[Qui73]{MR0338129}
D.~Quillen.
\newblock Higher algebraic {$K$}-theory. {I}.
\newblock In {\em Algebraic $K$-theory, I: Higher $K$-theories (Proc. Conf.,
  Battelle Memorial Inst., Seattle, Wash., 1972)}, pages 85--147. Lecture Notes
  in Math., Vol. 341. Springer, Berlin, 1973.

\bibitem[Qui78]{MR493916}
D.~Quillen.
\newblock Homotopy properties of the poset of nontrivial {$p$}-subgroups of a
  group.
\newblock {\em Adv. in Math.}, 28(2):101--128, 1978.

\bibitem[Seg68]{MR0232393}
G.~Segal.
\newblock Classifying spaces and spectral sequences.
\newblock {\em Inst. Hautes \'Etudes Sci. Publ. Math.}, (34):105--112, 1968.

\end{thebibliography}
\end{document}